\documentclass[11pt]{article}
\usepackage[a4paper,left=2.0cm,right=2.0cm,top=3cm,bottom=3cm]{geometry}
\usepackage{amsmath,amsfonts,amssymb}
\usepackage{mathtools}
\usepackage{siunitx}

\usepackage{fancyhdr}
\usepackage{graphicx}
\usepackage{float}
\usepackage{subfig}

\usepackage{booktabs}
\usepackage{multirow}

\def\@maketitle{\newpage
  \null
  \vskip 1em \begin{center} \let \footnote \thanks
    {\normalsize \bf \@title \par} \vskip 1em {\normalsize \lineskip .5em
	  \begin{tabular}[t]{c} \@author
      \end{tabular}\par}
	\vskip 1em {\normalsize \@date} \end{center}
	\par
  \vskip 1em}

\def\*#1{\mathbf{#1}}
\def\nsd{{n_\mathrm{sd}}}
\def\nts{{n_\mathrm{ts}}}

\def\nel{{n_\mathrm{el}}}

\def\ndf{{n_\mathrm{dof}}}

\newcommand{\bg}{{\boldsymbol{\mathit{g}}}}
\newcommand{\bh}{{\boldsymbol{\mathit{h}}}}

\newcommand{\bn}{\mathbf{n}}
\newcommand{\bp}{\mathbf{p}}

\newcommand{\bu}{\mathbf{u}}

\newcommand{\bx}{\mathbf{x}}

\newcommand{\bA}{\mathbf{A}}
\newcommand{\bB}{\mathbf{B}}

\newcommand{\bF}{\mathbf{F}}
\newcommand{\bG}{\mathbf{G}}

\newcommand{\bI}{\mathbf{I}}
\newcommand{\bJ}{\mathbf{J}}
\newcommand{\bK}{\mathbf{K}}
\newcommand{\bL}{\mathbf{L}}
\newcommand{\bM}{\mathbf{M}}

\newcommand{\bP}{\mathbf{P}}

\newcommand{\bS}{\mathbf{S}}

\newcommand{\bU}{\mathbf{U}}

\newcommand{\bW}{\mathbf{W}}

\newcommand{\bY}{\mathbf{Y}}

\newcommand{\zero}{\mathbf{0}}
\newcommand{\bxi}{{\boldsymbol{\xi}}}

\newcommand{\btau}{\boldsymbol{\tau}}
\newcommand{\ddx}[2]{\frac{\partial #1}{\partial #2}}
\newcommand{\ddt}[3]{\left(\ddx{#1}{#2}\right)^{#3}}
\newcommand{\tri}[1]{$\bigtriangleup_{#1}$}

\begin{document}

{\large \bf Simplex Space-Time Meshes in Compressible Flow Simulations} \par \vspace{6pt}
{\large Max von Danwitz\footnote{Chair for Computational
Analysis of Technical Systems (CATS),
RWTH Aachen University, 52056 Aachen,
Germany. Email: \{danwitz, karyofylli, hosters, behr\}@cats.rwth-aachen.de}, Violeta Karyofylli\footnotemark[1], Norbert Hosters\footnotemark[1], Marek Behr\footnotemark[1]} \par \vspace{6pt}

\begin{abstract}
\noindent
Employing simplex space-time meshes enlarges the scope of compressible flow simulations. The simultaneous discretization of space and time with simplex elements extends the flexibility of unstructured meshes from space to time. In this work, we adopt a finite element formulation for compressible flows to simplex space-time meshes. The method obtained allows, e.g., flow simulations on spatial domains that change topology with time. We demonstrate this with the two-dimensional simulation of compressible flow in a valve that fully closes and opens again. Furthermore, simplex space-time meshes facilitate local temporal refinement. A three-dimensional transient simulation of blow-by past piston rings is run in parallel on 120 cores. The timings point out savings of computation time gained from local temporal refinement in space-time meshes.
\end{abstract}

\section{Introduction}
\label{sec:intro}
In this work, we combine two constituents. A finite element formulation for compressible Navier--Stokes equations and simplex space-time meshes.
 
%Space-time finite element simulations of compressible flows have a long tradition and can be traced back to the publication of Jamet and Bonnerot in 1975~\cite{jamet1975numerical}.
Regarding the first constituent, our formulation can be traced back to the extension of the streamline upwind/Petrov-Galerkin (SUPG) method~\cite{brooks1982streamline} to compressible flow problems.\cite{hughes1984finite} The formulation based on conservation variables as primary unknowns posed a breakthrough in handling the advective character of compressible flows on equal-order interpolation elements. Building on this, Shakib in his thesis~\cite{shakib1989finite} developed and analyzed space-time finite element algorithms for compressible flow simulations using entropy variables. Entropy variables as primary unknowns lead to a symmetric from of the compressible flow equations. Guidance in the question of which variable set to use as primary unknowns can be found in the comparative study of Hauke.\cite{hauke1998comparative} The performance of four varialbe sets, namely conservation, entropy, density-, and pressure-primitive variables, has been evaluated for various test cases. Pressure-primitive variables were found to be well behaved in the incompressible limit and the most convenient to prescribe boundary conditions. Recently, Xu et al. further developed the formulation based on pressure-primitive variables and extended it for moving domains, weakly enforced essential boundary conditions, and sliding interfaces, with a generalized $\alpha$-method for time integration.\cite{xu2017compressible} In the following sections, we will adopt essential parts of this formulation for simplex space-time finite elements.

Simplex space-time meshes as the second constituent of our work provide a simultaneous discretrization of space and time. The meshes are characterized by an unstructured discretization of the spatial and temporal extent of the computational domain. In 1988, the potential of unstructured space-time meshes has been first postulated by Hughes and Hulbert.\cite{hughes1988space} They proposed solving a one-dimensional elastodynamics problem efficiently on an unstructured space-time mesh. Later, Maubach numerically solved the two-dimensional advection-diffusion equations on unstructured space-time meshes and provided a mathematical analysis of the method.\cite{maubach1991iterative} Unstructured space-time meshes were also applied in other fields, such as viscoelastic problems\cite{idesman2001finite} for example. However, the spatial domains of the transient problems remained two-dimensional for roughly 20 years.

The key step in using unstructured space-time meshes for transient three-dimensional problems is the mesh generation. A robust and simple procedure to construct four-dimensional simplex-based space-time meshes has been presented by Behr.\cite{behr2008simplex} The full space-time cylinder is first subdivided into space-time slabs which are in turn discretized with prism-type elements, followed by a Delaunay triangulation of each prism-type element. Neumüller and Steinbach proposed an alternative approach to four-dimensional simplex mesh generation, which includes the refinement of space-time meshes based on the Freudenthal algorithm.\cite{neumuller2011refinement} A third four-dimensional meshing strategy based on combinatorics was described in the thesis of Wang.\cite{wang2015discontinuous}

With mesh generation procedures available, recently three-dimensional transient problems were solved harnessing the advantages of unstructured space-time meshes. We describe here three examples in the field of fluid dynamics. A high-order discontinuous Galerkin method for compressible flows on domains with large deformations~\cite{wang2015high} was extended to three dimensions.\cite{wang2015discontinuous} Also in conjunction with a discontinuous Galerkin method, Lehrenfeld employed unstructured space-time discretizations to solve two-phase mass transport problems in three dimensions.\cite{Lehrenfeld2015On} As a third example, the incompressible two-phase flow simulations of Karyofylli et al. were boosted by adaptive temporal refinement along the moving fluid interface using simplex finite elements.\cite{karyofylli2018simplex} Also linear algebra aspects of unstructured space-time discretizations were recently explored. Steinbach and Yang studied the performance of algebraic multigrid preconditioned GMRES methods in solving linear systems arising from adaptive three- and four-dimensional space-time discretizations of the heat equation.\cite{steinbach2018comparison}

Now, with both constituents introduced, we proceed with the combination as follows. In Section~\ref{sec:cns}, we briefly review the Navier--Stokes equations of compressible flows and recast them as generalized advective-diffusive system. In Section~\ref{sec:numerics}, we discuss the numerical treatment and present a formulation suitable for simplex space-time finite elements. Section~\ref{sec:examples} aims at experimentally validating the formulation, as well as, demonstrating the particular potential of compressible flow simulations on simplex space-time meshes. Section~\ref{sec:conclusion} contains concluding remarks.

\section{Governing equations of compressible flows}
\label{sec:cns}
\subsection{Preliminaries}

As the Navier--Stokes equations of compressible flows state the pointwise conservation of mass, momentum, and energy, they are naturally expressed in terms of the conserved quantities, called conservation variables $\mathbf{\tilde{U}}$. For a three-dimensional problem the conservation varialbes are
\begin{equation}
\mathbf{\tilde{U}} = \left[ \rho, \rho u_1, \rho u_2, \rho u_3, \rho e_{tot}\right]^\top, 
\end{equation}
where $\rho$ is the density, $u_i$ are the velocity components in each space direction, $i = 1, ... , \nsd$, collected in $\bu$, and $e_{tot} = e + \frac{1}{2} \| \bu\|^2 $ is the total energy per unit mass of the fluid, being the sum of the internal energy $e$ and the kinetic energy $\frac{1}{2} \| \bu\|^2 $.

An alternative set of five independent variables are the pressure-primitive variables
\begin{equation}
\mathbf{Y} = \left[  p,  u_1,  u_2, u_3, T \right]^\top,
\end{equation}
where $p$ is the pressure and $T$ is the temperature. Throughout this paper we will consider an ideal, calorically perfect gas with
\begin{equation}
p = \rho R T, \quad \text{and} \quad  e  = \frac{R}{\gamma -1} T.
\end{equation}
$R$ is the specific gas constant, and $\gamma$ the ratio of specific heats.

\subsection{Strong form}
Expressed in conservation variables, the strong form of the compressible Navier--Stokes equations is given as
\begin{equation}
\mathbf{\tilde{U}}_{,t} + \mathbf{\tilde{F}}^{\text{adv}}_{i,i} - \mathbf{\tilde{F}}^{\text{diff}}_{i,i} -  \mathbf{\tilde{S}} = \zero.
\label{eq:cons}
\end{equation}
Therein, $\mathbf{\tilde{S}}$ is the vector of source terms. The vector of advective fluxes $\mathbf{\tilde{F}}^{\text{adv}}_{i} $ and the vector of diffusive fluxes $\mathbf{\tilde{F}}^{\text{diff}}_{i}$ are defined as
\begin{equation}
\mathbf{\tilde{F}}^{\text{adv}}_{i} = \left[  \begin{array}{c}
	     \rho u_i \\
	     \rho u_i u_1 + p \delta_{1i} \\
	    \rho u_i u_2 + p \delta_{2i} \\
	    \rho u_i u_3 + p \delta_{3i} \\
	     \rho u_i e_{tot} + p u_i
	     \end{array} \right], \qquad 
	     \mathbf{\tilde{F}}^{\text{diff}}_{i} = \left[  \begin{array}{c}
	     0 \\
	     \tau_{1i} \\
	     \tau_{2i} \\
	     \tau_{3i} \\
	     \tau_{ij} u_j - q_i  
	     \end{array} \right].
\end{equation}

Here, $\delta_{ij}$ is the Kronecker delta. The viscous stress tensor $\btau$ models the stress contribution caused by the strain rate; its entries $\tau_{ij}$ are defined as
\begin{equation}
\tau_{ij} = \mu \left( u_{i,j} + u_{j,i} \right) +  \lambda u_{k,k} \delta_{ij},  \quad i,j = 1, \dots , \nsd,
\end{equation}
where $\mu$ is the dynamic viscosity and $ \lambda = -\frac{2}{3} \mu$ the bulk viscosity. The entries $q_i$ of the heat flux vector are
\begin{equation}
q_i =  - \kappa T_{,i},
\end{equation}
where $\kappa$ is the coefficient of thermal conductivity.
In the following, we denote partial time derivatives with $(\cdot)_{,t}$, partial derivatives in spatial directions with $(\cdot)_{,i}$. The Einstein summation convention applies to repeated indices.

\subsection{Reduced form of the energy equation}

By algebraic manipulation of the equation system~\eqref{eq:cons}, the energy equation of the compressible Navier--Stokes equations can be cast in a reduced form,\cite{xu2017compressible} leading to the following system
\begin{equation}
\mathbf{U}_{,t} + \mathbf{F}^{\text{adv \textbackslash p}}_{i,i} + \mathbf{F}^{\text{p}}_{i,i} + \mathbf{F}^{\text{sp}} - \mathbf{F}^{\text{diff}}_{i,i} - \mathbf{S} = \zero.
\label{eq:red}
\end{equation}
The conservation variables with reduced energy equation,
\begin{equation}
\bU = \left[  \rho,  \rho u_1,   \rho u_2, \rho u_3, \rho e \right]^\top,
\end{equation}
differ from $\mathbf{\tilde{U}}$ in the fifth entry. $\rho e$ is solely the internal energy per unit mass instead of $\rho e_{tot}$ composed of internal and kinetic energy per unit mass.

%\begin{equation}
%\bU = \left[  \begin{array}{c}
%	     \rho  \\
%	     \rho u_1 \\
%	    \rho u_2 \\
%	    \rho u_3 \\
%	     \rho e
%	     \end{array} \right].
%\end{equation}

To facilitate a separate treatment in the weak form, the advective fluxes are split into the contribution to the advective fluxes without the pressure contribution $\mathbf{F}^{\text{adv \textbackslash p}}_{i}$ and the pressure contribution $\mathbf{F}^{\text{p}}_{i}$. Further, the balance laws with reduce energy equation contain the contribution of stress power to the energy equation $\mathbf{F}^{\text{sp}}$ and the diffusive fluxes $\mathbf{F}^{\text{diff}}_{i}$. The four flux vectors read
\begin{equation}
\mathbf{F}^{\text{adv \textbackslash p}}_{i} = \left[  \begin{array}{c}
	     \rho u_i \\
	     \rho u_i u_1  \\
	    \rho u_i u_2  \\
	    \rho u_i u_3  \\
	     \rho u_i e
	     \end{array} \right], \qquad 
	     \mathbf{F}^{\text{p}}_{i} = \left[  \begin{array}{c}
	     0 \\
	     p \delta_{1i} \\
	     p \delta_{2i} \\
	    p \delta_{3i} \\
	    0
	     \end{array} \right],\qquad
\mathbf{F}^{\text{sp}} = \left[  \begin{array}{c}
	     0 \\
	     0 \\
	     0 \\
	     0 \\
	     p u_{i,i} - \tau_{i,j} u_{j,i}
	     \end{array} \right], \qquad 
	     \mathbf{F}^{\text{diff}}_{i} = \left[  \begin{array}{c}
	     0 \\
	     \tau_{1i} \\
	     \tau_{2i} \\
	     \tau_{3i} \\
	    -q_i
	     \end{array} \right].
\end{equation}
$\bS$ can be used to model source terms such as body forces in the momentum equations. 

\subsection{Generalized advective-diffusive system}

A formulation with $\rho$ as variable is not well-defined in the incompressible limit.\cite{hauke1998comparative} Also, the specification of boundary conditions is cumbersome when using the conservation variables $\bU$ as primary unknowns. Hence, we consider the conservation variables $\bU$ as function of the pressure-primitive variables $\bY$, leading to the change of variables $\bU=\bU(\bY)$.\cite{hauke1998comparative}
For convenience, the following mappings are introduced
\begin{equation}
\bA_0 \coloneqq \ddx{\bU}{\bY}, \quad \bA^{\text{adv \textbackslash p}}_i \coloneqq \ddx{\mathbf{F}^{\text{adv \textbackslash p}}_{i}}{\bY}, \quad \bA^{\text{p}}_i \coloneqq \ddx{\mathbf{F}^{\text{p}}_{i}}{\bY}.
\end{equation}
Further, the matrices $\bA^{\text{sp}}_i$ and $\bK_{ij}$ are constructed, such that
\begin{equation}
\bA^{\text{sp}}_i \bY_{,i} = \bF^{\text{sp}}, \quad \bK_{ij} \bY_{,j} = \bF^{\text{diff}}_i.
\end{equation}
Explicit forms of the matrices are given in the appendix of the publication by Xu et al.~\cite{xu2017compressible}

Using the change of variables $\bU=\bU(\bY)$,  and the mappings above, Equation~\eqref{eq:red} can be written as generalized advective-diffusive system
\begin{equation}
\textbf{Res}(\bY) \coloneqq \bA_0 \bY_{,t} + \left( \bA^{\text{adv \textbackslash p}}_i + \bA^{\text{p}}_i + \bA^{\text{sp}}_i \right) \bY_{,i} - (\bK_{ij} \bY_{,j})_{,i} -\bS = \zero.
\label{eq:genad}
\end{equation}

\section{Numerical treatment with simplex space-time finite elements}
\label{sec:numerics}

\subsection{Simplex finite elements}
\label{ssec:simplices}
\begin{figure}[h]
\centering
\includegraphics[width=0.8\textwidth]{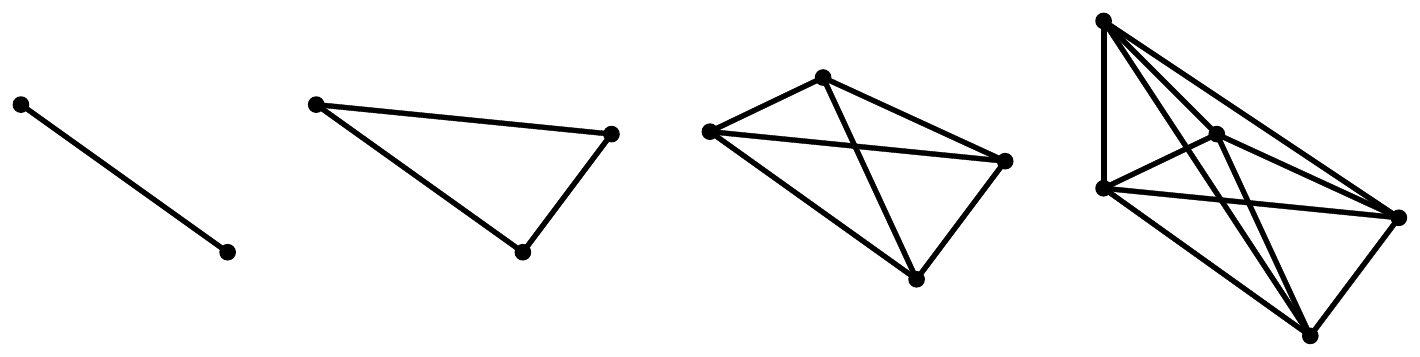}
\caption{Examples of Simplex Finite Elements for $d=1,2,3,$ and $4$ previously published in~\cite{karyofylli2018simplex}.}
\label{fig:simplices}
\end{figure}

A point $\bx \in \mathbb{R}^d$ is characterized by its coordinates $x_i, i =1,...,d$. As commonly done, we arrange the coordinates in a column vector
\begin{equation}
\bx = \left[ \begin{array}{c}
	x_1  \\
	\vdots \\
	x_d  \\
	\end{array} \right].
\end{equation}
In the case $d=2$, we denote $x_2$ by $y$. To discretize any finite-dimensional polytope region $R \subset \mathbb{R}^d$, without gaps, $d$-simplices can be used. Simplices for $d=1,2,3$, and $4$ are shown in Figure~\ref{fig:simplices}. The 1-simplex is a line, the 2-simplex a triangle, the 3-simplex a tetrahedron, and the 4-simplex a pentatope. Each $d$-simplex has $d+1$ vertices. We refer to the points where the simplex vertices of our $d$-dimensional tesselation meet as nodes.

%\begin{figure}[h]
%\centering
%\begin{tikzpicture}
%\node[anchor=west,scale=1.0] at (0,0) {\begin{tikzpicture} \input{../../CommonPics/PlainLine.tex} \end{tikzpicture}};
%\node[anchor=west,scale=1.0] at (3,0)  {\input{../../CommonPics/PlainTri.tex}};
%\node[anchor=west,scale=1.0] at (7,0)  {\input{../../CommonPics/PlainTet.tex}};
%\node[anchor=west,scale=1.0] at (11,0) {\input{../../CommonPics/PlainPent.tex}};
%\end{tikzpicture}
%\caption{Examples of simplices for $d=1,2,3,$ and $4$. Previously published by Karyofylli et al.~\cite{karyofylli2018simplex}}
%\label{fig:simplices}
%\end{figure}

\begin{figure}
\centering
\includegraphics[width=0.8\textwidth]{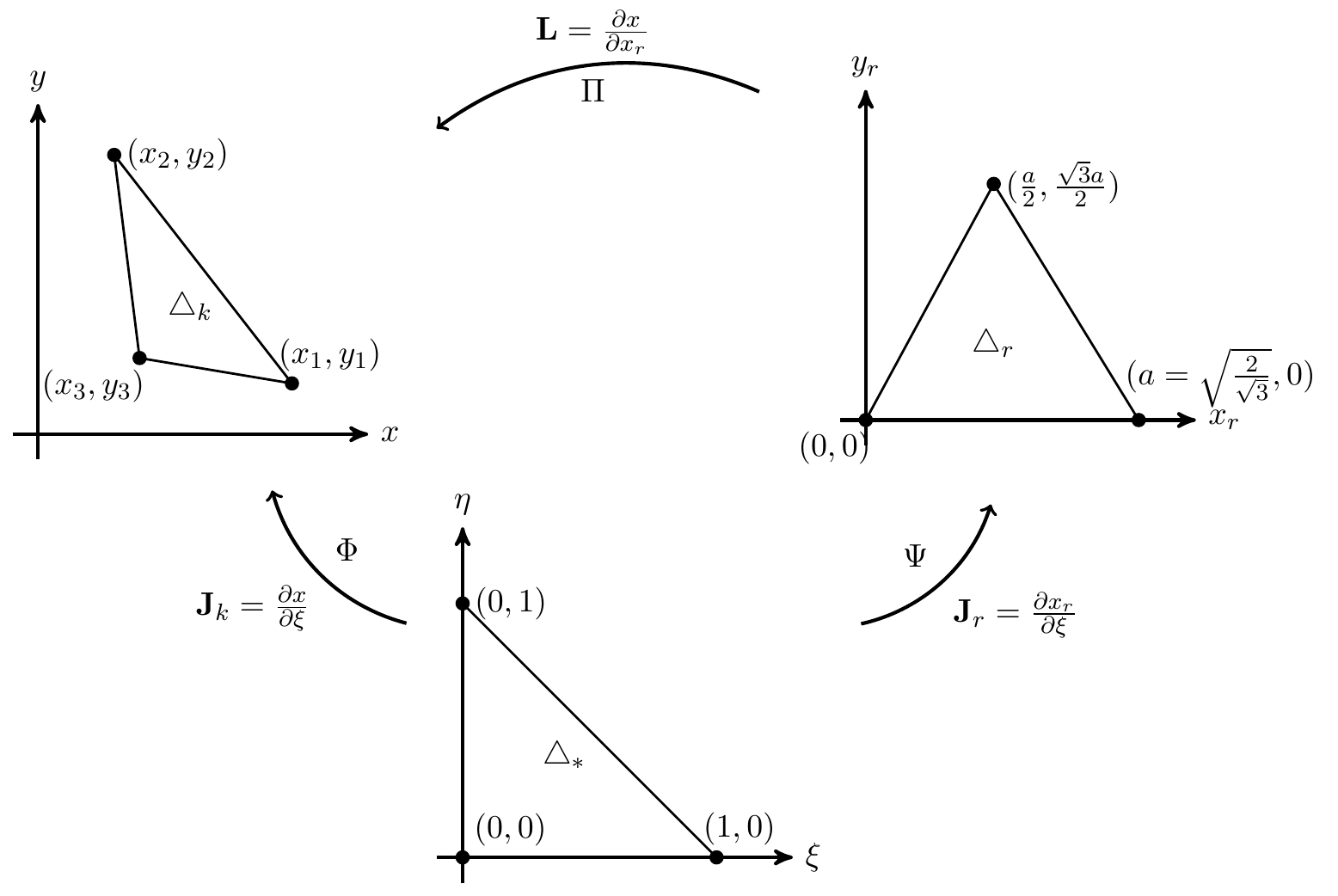}
\caption{Mappings between $\mathbb{P}_1$ reference element \tri{*}, physical element \tri{k} and regular element \tri{r} in the two-dimensional case.}
\label{fig:mappings}
\end{figure}

We further define the isoparametric local-global mapping $\Phi$ from a reference element \tri{*} with coordinates $\bxi$ to a physical element \tri{k} with coordinates $\bx$ as
\begin{equation}
\bx =\Phi(\bxi) = \sum_{q=1}^{d+1} \bx_q \chi_q(\bxi),
\end{equation}
where the index $q$ runs over the nodes of the physical element \tri{k}.\cite{elman2014finite} See Figure~\ref{fig:mappings} for a two-dimensional visualization. In two dimensions, the mapping reads
\begin{eqnarray}
x(\xi,\eta) &=& x_1 \chi_1(\xi,\eta) + x_2 \chi_2(\xi,\eta) + x_3 \chi_3(\xi,\eta),\nonumber \\
\label{eq:x}
y(\xi,\eta) &=& y_1 \chi_1(\xi,\eta) + y_2 \chi_2(\xi,\eta) + y_3 \chi_3(\xi,\eta),
\end{eqnarray}

with the $\mathbb{P}_1$ basis functions defined on the reference element \tri{*}
\begin{eqnarray}
\chi_1(\xi,\eta) &=& 1-\xi-\eta, \nonumber \\
\label{eq:P1}
\chi_2(\xi,\eta) &=& \xi, \\ \nonumber
\chi_3(\xi,\eta) &=& \eta \nonumber.
\end{eqnarray}

Generalized to $d$-dimensions, the $\mathbb{P}_1$ basis functions are
\begin{eqnarray}
\chi_1(\bxi) &=& 1 - \sum_{i=1}^{d}\xi_i \nonumber  \\
\chi_j(\bxi) &=& \xi_j,\quad  j=2, \dots, d+1.
\label{eq:P1d}
\end{eqnarray}

We follow the usual convention and define the Jacobian matrix of a vector-valued function by varying the differential operator, e.g., $ \ddx{}{\xi} $, between columns. In two dimensions, the Jacobian matrix reads

\renewcommand{\arraystretch}{1.2}
\begin{equation}
	\bJ = \left( \begin{array}{cc}
	\ddx{x}{\xi} & \ddx{x}{\eta}  \\
	\ddx{y}{\xi} & \ddx{y}{\eta}  \\
	\end{array} \right).
	\label{eq:J}
\end{equation}
Inserting the $\mathbb{P}_1$ basis function definitions~\eqref{eq:P1} into~\eqref{eq:x} and that into~\eqref{eq:J}, $\bJ_k$ can be computed for a given two-dimensional physical element \tri{k} as
\begin{equation}
	\bJ_k = \ddx{(x,y)}{(\xi,\eta)} = \left( \begin{array}{cc}
	x_2-x_1 & x_3-x_1  \\
	y_2-y_1 & y_3-y_1  \\
	\end{array} \right).
\end{equation}
This generalizes to $d$-dimensions as
\begin{equation}
	\bJ_k = \ddx{\bx}{\bxi} = \left( \begin{array}{cccc}
	\bx_2-\bx_1, & \bx_3-\bx_1,  & \dots,  & \bx_{d+1} - \bx_1\\
	\end{array} \right).
	\label{eq:Jd}
\end{equation}
The explicit form of the Jacobian confirms that $\bJ_k$ is constant per element, which shows a characteristic of the isoparametric $\mathbb{P}_1$ mapping $\Phi$ being an affine transformation. We further observe, that the columns of $\bJ_k$ are the distance vectors $\bp_i = \bx_{i+1} - \bx_1$ from node 1 to the other $d$ nodes of the simplex. Clearly, permuting the node numbering changes $\bJ_k$.

\subsection{Metric tensor definition}
\label{ssec:metric}
In stabilized finite element flow simulations, the covariant metric tensor $\hat{\bG}$ is frequently used to incorporate directional element length information into the stabilization parameter.\cite{shakib1989finite,hauke1998comparative,xu2017compressible,pauli2017stabilized} (We will discuss the details of stabilization in Section~\ref{ssec:supg}.) The definition of the covariant metric tensor $\hat{\bG}$ is based on the inverse Jacobian $\bJ_k^{-1} = \ddx{\bxi}{\bx}$
\begin{equation}
\hat{G}_{ij} = \sum_{k=1}^d \ddx{\xi_k}{x_i} \ddx{\xi_k}{x_j}, \quad i,j= 1, \dots , d.
\end{equation}

However, with this definition, $\hat{\bG}$ varies under permutations of the physical element's node numbering for simplex elements. To achieve invariance, we incorporate a regular simplex element \tri{r} in the definition of the metric tensor.\cite{pauli2017stabilized} The $d$-dimensional regular simplex element~\tri{r} is constructed, such that all its edges are of equal length $l$. Requiring further that \tri{r} has the volume of the reference element
\begin{equation}
\text{vol(\tri{*})} = \frac{1}{d!} \overset{!}{=} \frac{l^d}{d!} \sqrt{\frac{d+1}{2^d}} = \text{vol(\tri{r})}
\label{eq:vol}
\end{equation}
leads to the dimension-dependent edge length of
\begin{equation}
l=\sqrt{2(d+1)^{-\frac{1}{d}}}.
\label{eq:smalll}
\end{equation}
The formulas for the volumes in Equation~\eqref{eq:vol} can be proven by induction.\cite{buchholz1992perfect}

Utilizing the $\mathbb{P}_1$ basis functions defined on the reference element~\tri{*}, we define a mapping $\Psi$ from the reference element \tri{*} to the regular simplex element \tri{r}, i.e., $\bx_r = \Psi(\bxi)$. The Jacobian $\bJ_r = \ddx{\bx_r}{\bxi}$ associated with this mapping can be computed in the same way as $\bJ_k$ using~\eqref{eq:Jd}. From vol(\tri{r}) $ \overset{!}{=}$ vol(\tri{*}), follows det$(\bJ_r) = 1$. The composition of $\Phi$ and $\Psi^{-1}$ we call $\Pi = \Phi \circ \Psi^{-1}$, which maps from the regular simplex element \tri{r} to the physical element \tri{k}. A visualization can be found in Figure~\ref{fig:mappings}. The Jacobian $\bL$ associated with the mapping $\Pi$ is
\begin{equation}
\bL =  \ddx{\bx}{\bx_r} = \ddx{\bx}{\bxi} \ddx{\bxi}{\bx_r}.
\end{equation}

If we base the definition of the contravariant metric tensor $\bG^{-1}$ on $\bL$, then
\begin{equation}
\bG^{-1} = \bL \bL^{\top} = \ddx{\bx}{\bx_r} \left(\ddx{\bx}{\bx_r}\right) ^\top = \ddx{\bx}{\bxi} \ddx{\bxi}{\bx_r} \ddt{\bxi}{\bx_r}{\top} \ddt{\bx}{\bxi}{\top},
\label{eq:contra}
\end{equation}
is invariant under cyclic permutations of the physical element's node numbering, and so is the covariant metric tensor $\bG$ defined as its inverse
\begin{equation}
\bG = \left(\bG^{-1}\right)^{-1} = \ddt{\bx}{\bxi}{-\top}\ddt{\bxi}{\bx_r}{-\top} \ddt{\bxi}{\bx_r}{-1} \ddt{\bx}{\bxi}{-1} = \ddt{\bxi}{\bx}{\top} \ddt{\bx_r}{\bxi}{\top} \ddx{\bx_r}{\bxi} \ddx{\bxi}{\bx}.
\end{equation}

The mapping $\ddx{\bx_r}{\bxi}$ is known for each $d$-simplex and can be precomputed, such that the computation of $\bG$ simplifies to
\begin{equation}
\bG = \ddt{\bxi}{\bx}{\top} \bM \,  \ddx{\bxi}{\bx}, \quad \text{with} \quad \bM = \ddt{\bx_r}{\bxi}{\top} \, \ddx{\bx_r}{\bxi} = \bJ_r^\top \bJ_r.
\label{eq:cov}
\end{equation}

$ \bM = \bJ_r^\top \bJ_r$ prompts the observation that the entries of $\bM$ are in fact scalar products of the distance vectors between the nodes of the regular simplex (compare Equation~\eqref{eq:Jd})
\begin{equation}
M_{ij} = \bp_i \cdot \bp_j = | \bp_i | |\bp_j | \cos(\alpha_{ij}).
\end{equation}
By construction, all edges of the regular simplex are of length $l$ (given in Equation~\eqref{eq:smalll}) and all angles between two edges are 60 degree. Hence,
\begin{equation}
|\bp_i| = l \quad \text{and} \quad \cos(\alpha_{ij}) = \frac{1}{2}+\frac{1}{2}\delta_{ij}
\end{equation}
hold for all choices of the base node for the distance vectors, i.e., $\bx_1$ in Equation~\eqref{eq:Jd}. It follows that $\bM$, characterized by
\begin{equation}
M_{ij} = l^2 (\frac{1}{2}+\frac{1}{2}\delta_{ij}),
\label{eq:M}
\end{equation}
exclusively measures edge lengths and angles of the regular simplex and therefore is invariant under permutations of the node numbering.
For $d=2,3,4$, $\bM$ reads, respectively
\begin{equation}
	\bM_{d=2} = \frac{1}{\sqrt{3}} \left( \begin{array}{cc}
	2 & 1  \\
	1 & 2  \\
	\end{array} \right), \quad%
		\bM_{d=3} = \frac{1}{\sqrt[3]{4}} \left( \begin{array}{ccc}
	2 & 1 & 1  \\
	1 & 2 & 1 \\
	1 & 1 & 2 \\
	\end{array} \right), \quad%
			\bM_{d=4} = \frac{1}{\sqrt[4]{5}} \left( \begin{array}{cccc}
	2 & 1 & 1  & 1\\
	1 & 2 & 1 & 1\\
	1 & 1 & 2 & 1\\
	1 & 1 & 1 &2 \\
	\end{array} \right).
\end{equation}

The explanation why $\bG$ is invariant under node permutations is now straightforward, we first observe that a permutation of the node numbering of the physical element is equivalent to a permutation of the node numbering of the reference element or the regular simplex element. Applying the node permutation to regular simplex element one obtains
\begin{equation}
\mathbf{\tilde{L}} = \ddx{\bx}{\bx_{\tilde{r}}} = \ddx{\bx}{\bxi} \ddx{\bxi}{\bx_{\tilde{r}}},
\end{equation}
which is associated with the mapping from the regular simplex element with permuted node numbering to the physical element. The corresponding covariant metric tensor with the permutation applied reads
\begin{equation}
\mathbf{\tilde{G}} = \left(\mathbf{\tilde{L}}\mathbf{\tilde{L}}^\top\right)^{-1} = \ddt{\bxi}{\bx}{\top} \ddt{\bx_
{\tilde{r}}}{\bxi}{\top} \ddx{\bx_{\tilde{r}}}{\bxi} \ddx{\bxi}{\bx} = \ddt{\bxi}{\bx}{\top} \mathbf{\tilde{M}} \,  \ddx{\bxi}{\bx} = \ddt{\bxi}{\bx}{\top} \bM \,  \ddx{\bxi}{\bx}.
\end{equation}
The last equality follows from the characterization of $\bM$ in Equation~\eqref{eq:M} and recovers the original definition of $\bG$ in Equation~\eqref{eq:cov}, showing that $\bG$ is indeed invariant under permutations of the physical element's node numbering. We will use the invariant metric tensor $\bG$ to construct a stabilization matrix $\btau$ that is likewise invariant under permutations of the element's node numbering (see Section~\ref{ssec:supg}).

\subsection{Space-time discretization methods}
\label{ssec:st}
\begin{figure}
\centering
\subfloat[Flat space-time\label{fig:STFST}]{\includegraphics{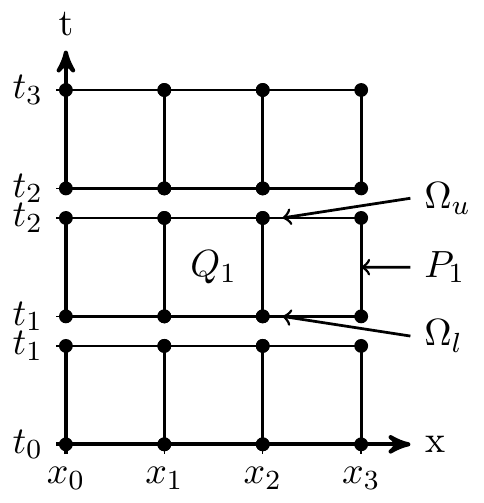}}\,
\subfloat[Simplex space-time\label{fig:STSST}]{\includegraphics{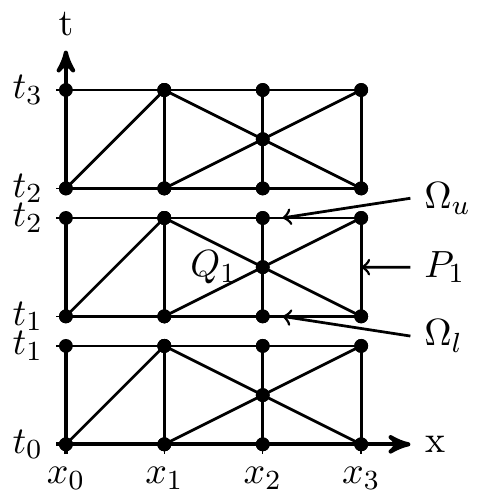}}\,
\subfloat[Unstructured space-time\label{fig:STUST}]{\includegraphics{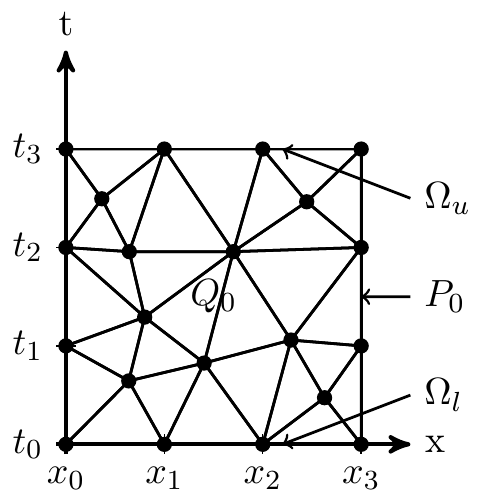}}\\
\caption{Space-time discretization methods.}
\label{fig:ST}
\end{figure}

With the mappings and metric tensors for simplex elements defined, we can now proceed to discretize the subset of the space-time continuum $Q$ in which we wish to solve the compressible Navier--Stokes equations. Thinking of $Q$ being generated by the time dependent spatial domain $\Omega_t \subset \mathbb{R}^{\nsd}$ evolving over the time interval $I \subset \mathbb{R}$ shows that $Q \subset \mathbb{R}^{\nsd+1}$.  Based on the discontinuous Galerkin method in time, $Q$ is sliced into $N$ space-time slabs $Q_n$. The extent of the space-time slabs in temporal direction is denoted by $\Delta t$. An exemplary slicing of $Q =[x_0, x_3]\times ]t_0,t_3]$ is shown in Figure~\ref{fig:STFST}. Each space-time slab is bounded by the spatial domain at the lower and upper time level $\Omega_{\text{l}}$, $\Omega_{\text{u}}$, respectively, as well as by $P \subset \mathbb{R}^{\nsd}$, which is the temporal evolution of the spatial domain boundary $\Gamma \subset \mathbb{R}^{\nsd-1}$.

In this paper, three methods to discretize $Q_n$ are considered. The flat space-time (FST) method extrudes a spatial discretization of $\Omega$ in time generating a discretization of $Q_n$ with prismatic space-time elements $Q^e_n$ (Figure~\ref{fig:STFST}). In the simplex space-time (SST) method, these prismatic elements are further subdivided into simplex elements $Q^e_n$ (Figure~\ref{fig:STSST}), allowing for local temporal refinement.\cite{behr2008simplex} A third option is the unstructured space-time (UST) method, which generates an unstructured tesselation of $Q_n$ with simplex elements $Q^e_n$, e.g., by Delaunay triangulation (Figure~\ref{fig:STUST}).

In Figure~\ref{fig:ST}, the spatial domain is restricted to one dimension for the sake of a clear visualization. However, the FST and SST discretizations are well established for complex three-dimensional spatial geometries, i.e., $Q \subset  \mathbb{R}^{4}$.\cite{karyofylli2018simplex} In case of UST, simplex discretizations of arbitrary $Q \subset \mathbb{R}^{3}$ can be obtained with a variety of commercially or freely available meshing tools. UST discretizations of $Q \subset \mathbb{R}^{4}$ of domains of engineering interest are to the authors' knowledge an open research problem.

\subsection{Weak form}
\label{ssec:weak}
Based on one of the above discretizations of the space-time slabs $Q_n$, trial and test functions are generated in each element $Q_n^e$ by the linear interpolation of the basis functions in the reference element and the mappings $(\Phi_n^e)^{-1}$ from physical elements to the reference element. Further, these functions are constructed to be $C^0$-continuous on the closure of the space-time slab $\overline{Q_n}$. Hence, we define the $H^1$-conformal finite element approximation space as
\begin{equation} 
H^1_{h,n} = \left\lbrace f^h \in  H^1(Q_n),\, f^h|_{Q^e_n} = v \circ (\Phi_n^e)^{-1} \; \text{with} \; v \in  \mathbb{V}_1,\, \forall Q_n^e \right\rbrace.
\end{equation}

In case of an FST discretization, $\mathbb{V}_1$ refers to the $\mathbb{PR}_1$ basis functions of the prismatic Lagrange finite element. In case of an SST or UST discretization, $\mathbb{V}_1$ refers to the $\mathbb{P}_1$ basis functions (Equation~\eqref{eq:P1d}) of the simplical Lagrange finite element.

Further, we define a projection $\bP$ that selects from $\bY^h \in [H^1_{h,n}]^{\nsd+2}$ the degrees of freedom, where Dirichlet boundary data $\bg^h$ is prescribed. The parts of the boundary $P_n$, where Dirichlet boundary conditions are prescribed, may vary for each of the primitive variables $\left( p,  u_1,  u_2, u_3, T \right)$. $\bP$ allows the definition of the trial function space
\begin{equation}
S_{h,n} = \left\lbrace \bY^h \in  [H^1_{h,n}]^{\nsd+2}, \bP \bY^h = \bg^h \right\rbrace 
\end{equation}
that fulfills the Dirichlet boundary conditions and the test function space

\begin{equation}
V_{h,n} = \left\lbrace \bW^h \in [H^1_{h,n}]^{\nsd+2}, \bP \bY^h = \zero  \right\rbrace
\end{equation}
with functions that vanish where Dirichlet boundary conditions are prescribed.

Considering that the finite element functions are discontinuous at the space-time slab boundaries $\Omega_l$ and $\Omega_u$, let $(\bY^h)^{\pm}_n$ abbreviate $\lim_{\varepsilon \to 0} \bY^h(t_n \pm \varepsilon)$.  Hence, the weak form of Equation~\eqref{eq:genad} can be stated as: 
For given initial conditions $(\bY^h)^-_0 = (\bY^h)_0$, find $\bY^h \in S_{h,n}$, such that on each time slab $Q_n$ and for all $\bW^h \in V_{h,n} $
\begin{equation}
B(\bW^h,\bY^h) = F(\bW^h).
\label{eq:Galerkin}
\end{equation}
Therein, the bilinear form reads
\begin{align}
		\label{eq:weak}
B(\bW^h,\bY^h) &=	\int_{Q_n} \bW^h \cdot \left[  \bA_0 \bY^h_{,t}  +  \left( \bA^{\text{adv \textbackslash p}}_i + \bA^{\text{sp}}_i \right) \bY^h_{,i} \right]  dQ \\ \nonumber
	&+	\int_{Q_n} \bW^h_{,i} \cdot \left[ \bK_{ij} \bY^h_{,j} - \bA_i^{\text{p}} \, \bY^h \right] dQ\\ \nonumber
	& + \int_{\Omega_{\text{l}}} \bW^h \cdot
		 \bA_0 \left( (\bY^h)^+_n - (\bY^h)^-_n \right)
		d\Omega ,
\end{align}
and the linear functional reads
\begin{align}
F(\bW^h) = \int_{Q_n} \bW^h \cdot \bS^h \, dQ +  \int_{P_n} \bW^h \cdot \bh^h\, dP.
\end{align}
Here we used
\begin{equation}
	     \bh^h = \left[  \begin{array}{c}
	     0 \\
	     -p n_1 +\tau_{1i} n_i\\
	     -p n_2 +\tau_{2i} n_i \\
	     -p n_3 + \tau_{3i} n_i \\
	    -q_i n_i
	     \end{array} \right],
\end{equation}
where $n_i$ are the components of the outwards pointing surface normal $\bn$, not to be confused with plain $n$ that indexes the space-time slabs form 0 to $N-1$.  For test cases with pentatope discretization, the domain boundary consists of tetrahedra in four dimensions. The normals to the boundary tetrahedra are computed using the generalized cross product definition.\cite{Lehrenfeld2015On} The diffusive and pressure fluxes in the second integral of Equation~\eqref{eq:weak} are integrated by parts and give rise to the boundary integral over $P_n$, while $\bh^h$ is evaluated with the solution of the previous Newton step. On the boundary between subsequent time slabs, i.e., $\Omega_l$, the third integral in Equation~\eqref{eq:weak} weakly enforces the coupling between the space-time slabs.

\subsection{SUPG operator}
\label{ssec:supg}
It is well-known that the Galerkin formulation for the convection-dominated compressible Navier--Stokes equations (Equation~\eqref{eq:Galerkin}) suffers from instabilities. To counteract these instabilities, the weak form is augmented by the following SUPG operator
\begin{equation}
B(\bW^h,\bY^h) +\bB_{\text{SUPG}}(\bW^h,\bY^h) = F(\bW^h),
\label{eq:nonlinear}
\end{equation}

\begin{equation}
\bB_{\text{SUPG}}(\bW^h,\bY^h) = \sum_{e=1}^{(\nel)_n} \int_{Q_n^e} \! 
	\left( (\hat{\bA}_{m})^{T} \bW^h_{,m}
	\right) \cdot \btau^e \mathbf{Res(Y}^h) dQ,
	\label{eq:supg}
\end{equation}
where $m =1, \dots, n_{sd}+1$ runs over $n_{sd}$ spatial dimensions and time. $\hat{\bA}_i$ are the combined advection matrices for conservation variables
\begin{equation}
\hat{\bA}_i = \left( \bA^{\text{adv \textbackslash p}}_i + \bA^{\text{p}}_i + \bA^{\text{sp}}_i \right)\bA_0^{-1}, \quad i = 1, \dots, \nsd.
\end{equation}
The advection matrix for conservation variables in time direction is the identity
\begin{equation}
\hat{\bA}_{t} = \bI.
\end{equation}
For the element-wise stabilization matrix $\btau^e$, several definitions are proposed in literature.\cite{shakib1989finite,xu2017compressible,hauke2001simple} We select the direct approach~\cite{xu2017compressible} and adopt it for unstructured space-time discretizations
\begin{align}
\btau^e = \left(\bG_{mp} \hat{\bA}_m \hat{\bA}_ p + C_{inv}^2 \bG_{ij} \bG_{kl} \hat{\bK}_{ik} \hat{\bK}_{lj} \right)^{-\frac{1}{2}}, \quad i,j,k,l = 1, \dots, \nsd, \quad m,p =1, \dots, \nsd+1.
\label{eq:tau}
\end{align}
The diffusion matrices for the conservation variables are obtained as $\hat{\bK}_{ij} = \bK_{ij} \bA_0^{-1} $. Furthermore, $\btau^e$ relies on the covariant metric tensor $\bG$ (constructed in Section~\ref{ssec:metric}). Using $\bG$, which is invariant under permutations of the element nodes, ensures that $\btau^e$ is invariant as well. In the FST case on non-moving domains, the full metric tensor $\bG$ can be decomposed into a spatial part $\left[ \bG \right]_{\nsd \times \nsd}$ and a temporal part $\left(\ddx{\theta}{t}\right)^2 = \frac{4}{\Delta t^2}$ as
\begin{equation}
\bG =  \left( \begin{array}{cccc}
	 &  &   & 0\\
    & \left[ \bG \right]_{\nsd \times \nsd} &  & \vdots \\
    &  &  & 0 \\
	0  &  \hdots & 0 &\left(\ddx{\theta}{t}\right)^2 \\
	\end{array} \right).
\end{equation}
For SST and UST discretizations, the full metric tensor $\bG$ is computed directly from the simplex space-time element with Equation~\eqref{eq:cov}.
The constant $C_{inv}$ scales the diffusive contribution to $\btau^e$. The choice of
\begin{equation}
C_{inv} = (\nsd+1)^2(\nsd+2)
\end{equation}
is motivated by an inverse estimate inequality for $\mathbb{P}_1$ simplex finite elements.\cite{knechtges2018simulation} The reference element considered therein is the simplex reference element of the spatial discretization common for FST, SST, and UST discretizations. In case of steady simulations, the indices $m$ and $p$ in Equations~\eqref{eq:supg} and ~\eqref{eq:tau} run only over the spatial dimensions. 

The principal square root of the $4\times4$ or  $5\times5$ matrix in Equation~\eqref{eq:tau} is computed using the direct Schur decomposition method~\cite{higham1987computing} provided through the NAG library.\cite{nagLibrary} The nonlinear equation system~\eqref{eq:nonlinear} is linearized using the Newton--Raphson technique. The resulting linear equation system is solved with a parallel, preconditioned GMRES implementation.\cite{behr1994finite}

\section{Numerical examples}
\label{sec:examples}
In the following numerical examples we will use the specific gas constant $R = 287 \frac{J}{kg K}$ and the ratio of specific heats $\gamma =1.4$. These are the values for air at standard conditions.\cite{anderson2010fundamentals} In the viscous computations in Sections~\ref{ssec:pulse},~\ref{ssec:valve}, and~\ref{ssec:pistonRings}, the temperature dependence of the dynamic viscosity $\mu(T)$ is modeled using Sutherland's relation
\begin{equation}
\mu = \mu_{ref} \frac{T_{ref}+C}{T+C}\left(\frac{T}{T_{ref}}\right)^{\frac{3}{2}},
\label{eq:Sutherland}
\end{equation}
with $\mu_{ref} = 21.7\cdot10^{-6}$ Pa s, $T_{ref}=  373.15$ K, $C=120$ K. The coefficient of thermal conductivity $\kappa$ is coupled to the dynamic viscosity $\mu$ using the Prandtl number $Pr = \frac{\mu}{\kappa}\frac{\gamma R}{\gamma-1} = 0.71 $.  For air and temperatures up to the order of 1000 K, it is sufficient to treat the Prandtl number as constant.\cite{anderson2010fundamentals}

\subsection{Supersonic flow over flat plate}
\label{ssec:plate}
\begin{figure}
\subfloat[Pressure isocontours]{\includegraphics[width=0.49\textwidth,trim={0cm 0cm 0cm 0cm},clip]{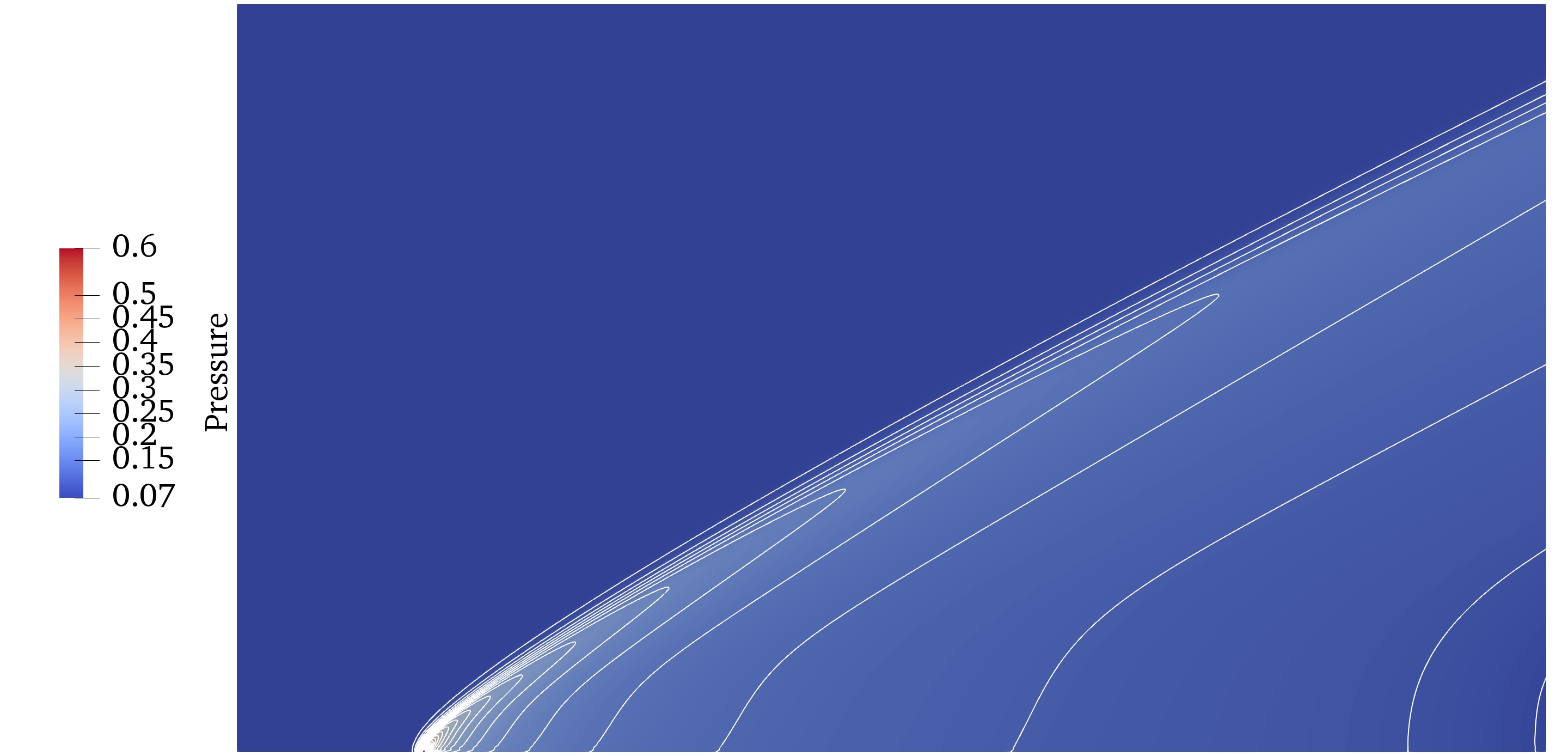}}
\subfloat[Velocity magnitude isocontours]{\includegraphics[width=0.49\textwidth,trim={0cm 0cm 0cm 0cm},clip]{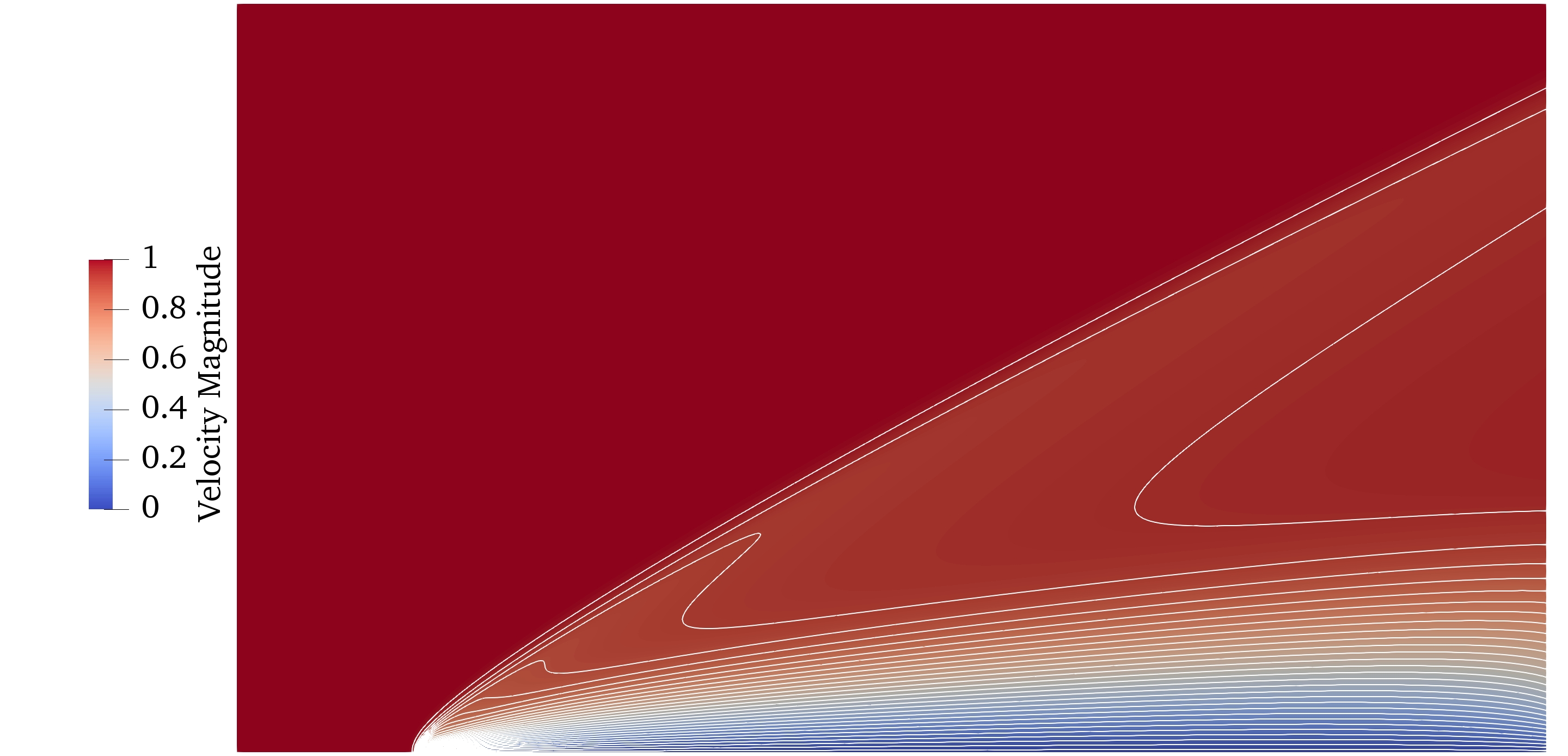}}\\
\subfloat[Temperature isocontours]{\includegraphics[width=0.49\textwidth,trim={0cm 0cm 0cm 0cm},clip]{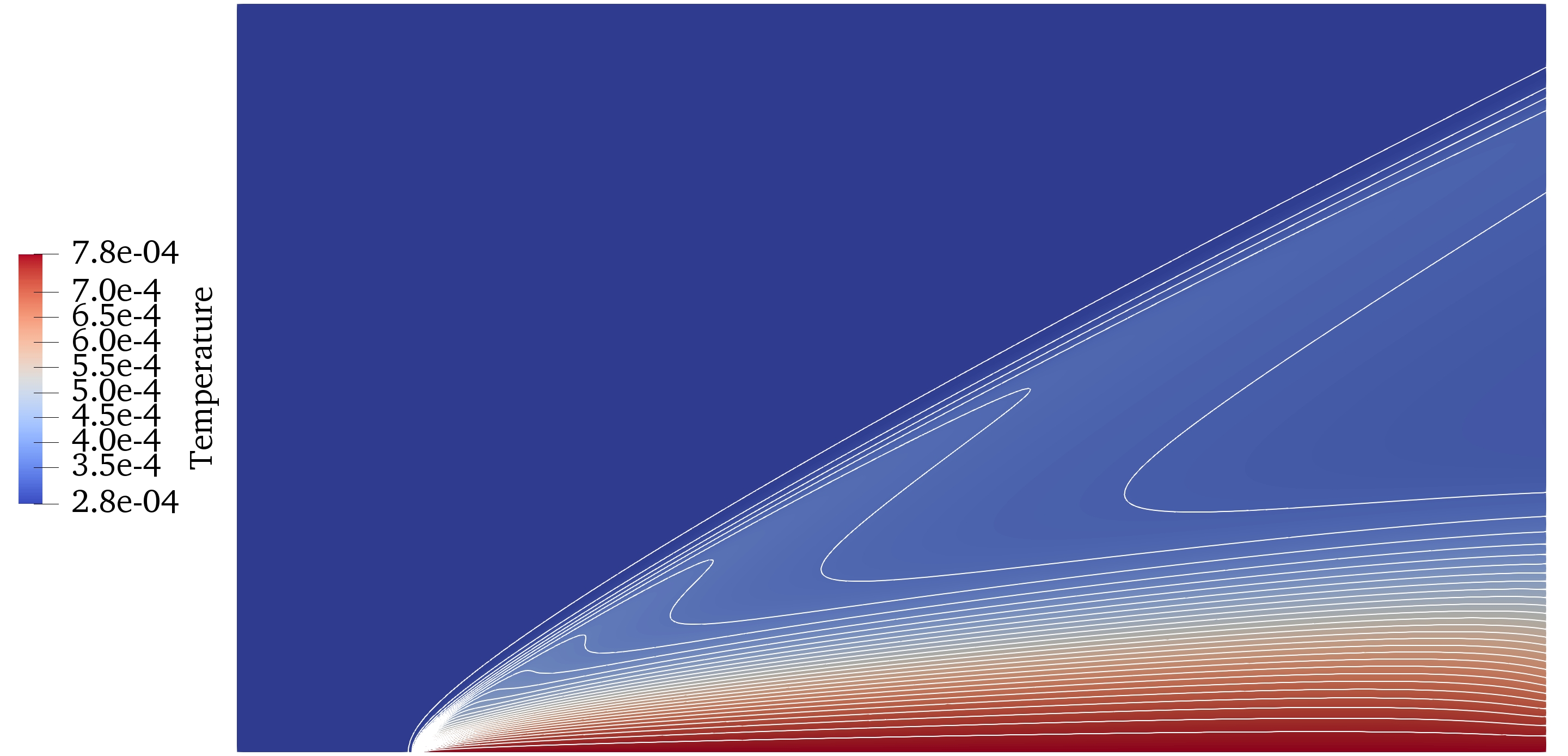}} \quad
\subfloat[Pressure coefficient along the solid wall \label{fig:fpPressure}]{\includegraphics[width=0.49\textwidth,trim={0cm 0cm 0cm 0cm},clip]{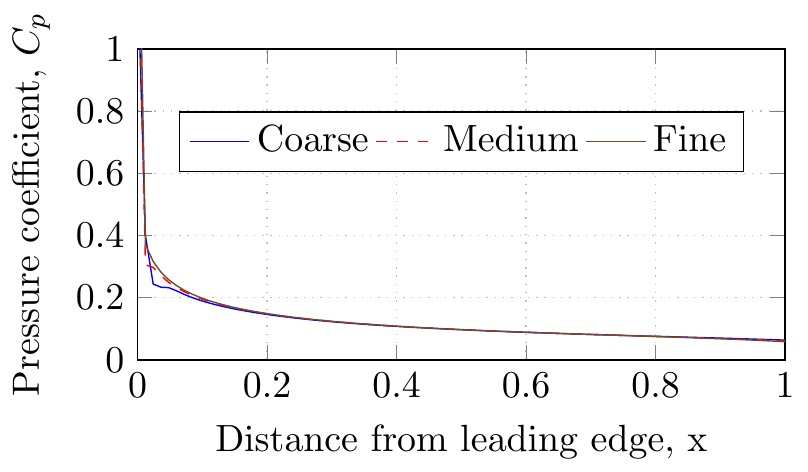}}
\caption{Supersonic flow over a flat plate. Flow field solutions and pressure coefficient.}
\label{fig:fpPics}
\end{figure}

The two-dimensional supersonic viscous flow over a flat plate at free-stream $Re_{\infty}=1000$ and $Ma_{\infty}= 3$ constitutes a standard test case for compressible flow simulations. The flow develops a curved shock and boundary layer along the flat plate (Figure~\ref{fig:fpPics}). Already in 1972 the test case was solved numerically.\cite{carter1972numerical} Later, the test case was used to study the stability and spatial accuracy of a space-time finite element formulation for compressilbe flows.\cite{shakib1989finite} More recently, it has been calculated with the following specifications:\cite{hauke1998comparative,xu2017compressible}

\begin{figure}
\centering
\subfloat[$\tau, q = 0$\label{fig:fpTa}]{\includegraphics[width=0.35\textwidth,trim={0cm 0cm 8cm 0cm},clip]{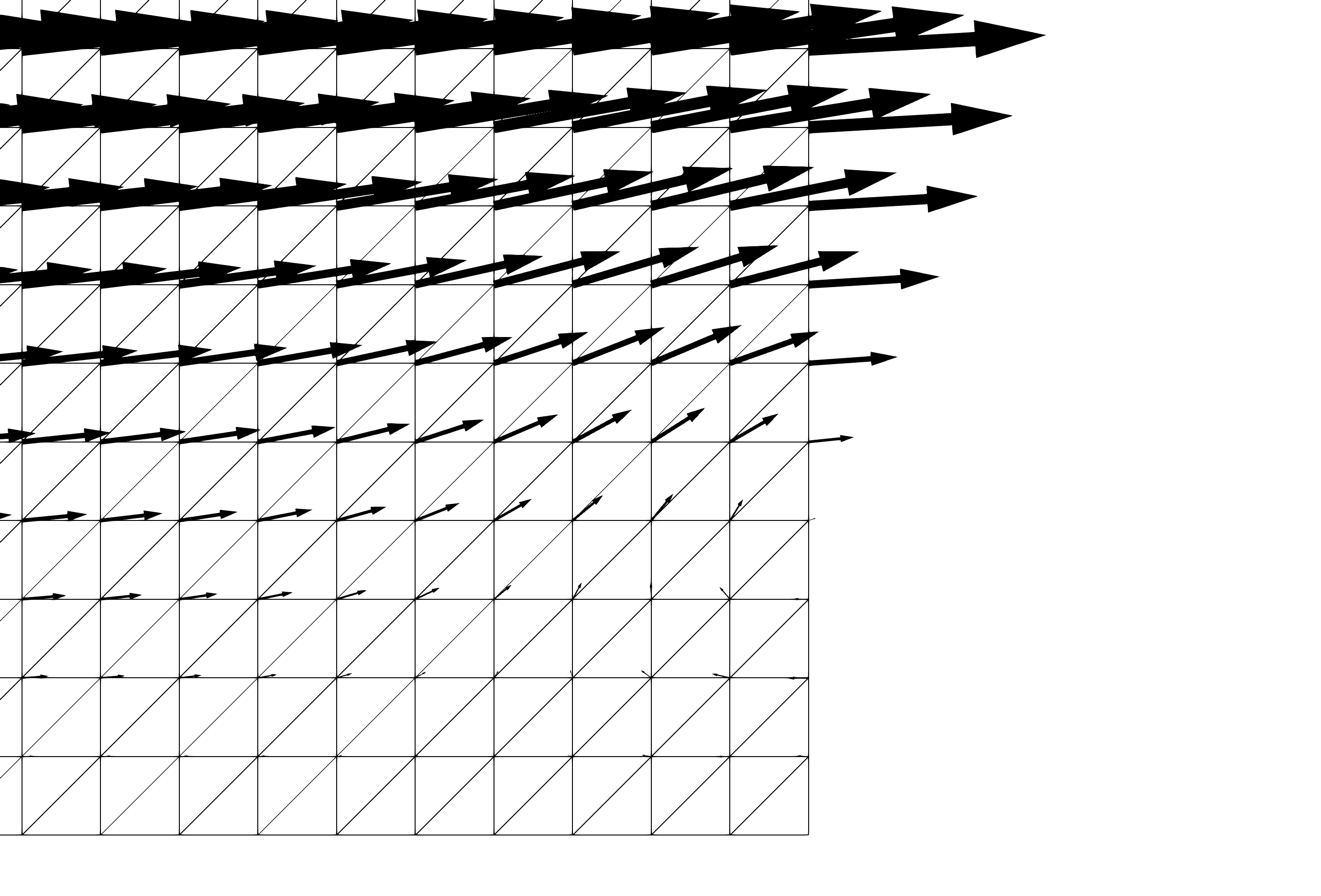}} \quad
\subfloat[$\tau, q \neq 0$\label{fig:fpTb}]{\includegraphics[width=0.35\textwidth,trim={0cm 0cm 8cm 0cm},clip]{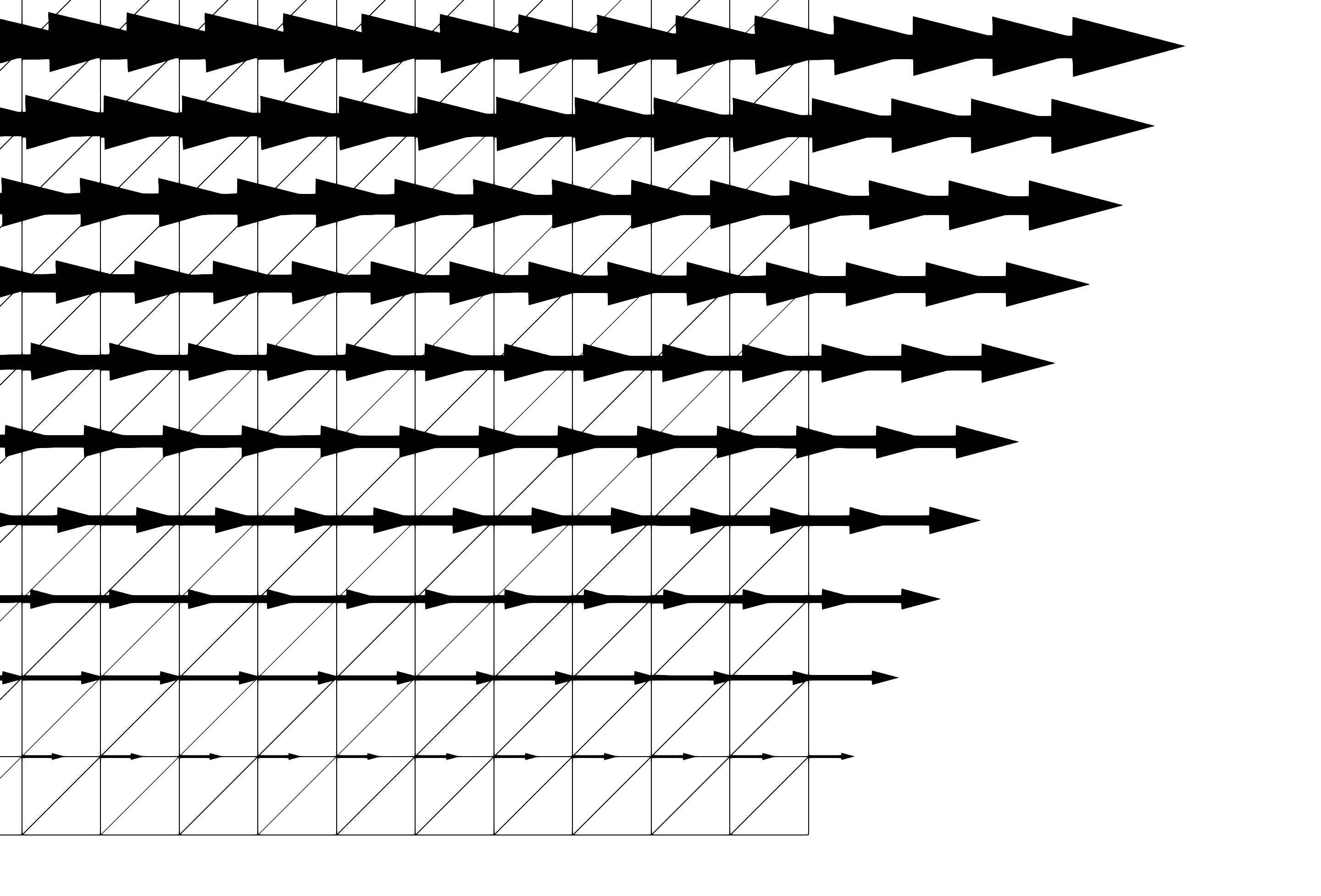}}\\
\caption{Supersonic flow over flat plate. Influence of the viscous terms in the boundary integral $\int_{P_n} \bW^h \cdot \bh^h \, dP$ along the outflow boundary on the velocity vectors next to the wall.}
\label{fig:fpT}
\end{figure}

The compressible Navier-Stokes equations are solved on a rectangle with $-0.2 \leq x \leq 1.2$ and $0 \leq y \leq 0.8$. A modified Sutherland's law
\begin{equation}
\mu = \frac{0.0906 \,T^{1.5}}{T+0.0001406},
\end{equation}
models the temperature dependence of the viscosity. On the left ($x=0.2$) and top ($y=0.8$) boundary, the free-stream conditions are prescribed for all degrees of freedom
\begin{equation}
	     \bY_{in} = \left[  \begin{array}{c}
	     p_{\infty} \\
	     u_1\\
	     u_2\\
	     T_{\infty} \\
	      \end{array} \right] = \left[  \begin{array}{c}
	     7.937 \times 10^{-2} \\
	     1 \\
	     0 \\
	     2.769 \times 10^{-4} \\
	      \end{array} \right].
\end{equation}
On the symmetry line ($y=0$ and $x<0$), the symmetry conditions $u_2=\tau_{12}=q_2=0$ are applied, rendering the viscous terms in the Neumann boundary integral $\int_{P_n} \bW^h \cdot \bh^h \, dP$ null. On the solid wall ($y=0$ and $x\geq 0$), the no-slip condition $\bu=\zero$ and the stagnation temperature of $T_W=T_{\infty}\left(1+\frac{\gamma -1}{2} {Ma_{\infty}}^2 \right) = 7.754\times10^{-4}$ are prescribed as Dirichlet boundary condition. Finally, along the outflow boundary on the right ($x=1.2$) no Dirichlet boundary condition is prescribed. However, we do include the viscous terms ($\tau, q \neq 0$) in the boundary integral $\int_{P_n} \bW^h \cdot \bh^h \, dP$ along the outflow boundary, which we find to stabilize the flow in the subsonic part of the outflow. Figure~\ref{fig:fpT} shows the velocity vectors in the bottom right corner of the computational domain, where the subsonic boundary layer touches the outflow boundary. Neglecting the viscous contributions to the boundary integral, we observe a small flow recirculation (see Figure~\ref{fig:fpTa}). Including the viscous contributions, a flow field with strictly aligned velocity vectors is obtained (see Figure~\ref{fig:fpTb}). Similar findings were reported for the compressible flow simulation with entropy variables.\cite{shakib1989finite}

%\begin{figure}
%\includegraphics[width=0.9\textwidth,trim={0cm 18cm 0cm 18cm},clip]{./FlatPlate/fpPressure}
%\includegraphics[width=0.9\textwidth,trim={0cm 18cm 0cm 18cm},clip]{./FlatPlate/fpVelocity}
%\includegraphics[width=0.9\textwidth,trim={0cm 18cm 0cm 18cm},clip]{./FlatPlate/fpTemperature}
%\caption{Supersonic flow over a flat plate. Isoconturs of pressure, velocity magnitude, and temperature.}
%\label{fig:fpPics}
%\end{figure}

To compare the results with the findings of Xu et al.,\cite{xu2017compressible} the computations are performed on three uniform meshes with 22400, 89600, and 358400 triangular elements, named 'Coarse', 'Medium' and 'Fine'. For the temporal discretization, the flat space-time method (FST) is employed. Figure~\ref{fig:fpPics} shows the isocontours of pressure, velocity magnitude, and temperature on the fine mesh. The curved shock and boundary layer are accurately resolved in all solution variables. The pressure coefficient $C_p = \frac{2(p-p_{\infty})}{\rho_{\infty}\|\bu_{\infty}\|^2}$ along the solid wall is shown for the three meshes in Figure~\ref{fig:fpPressure}. Except for a small kink close to the singularity at the leading edge of the plate, which vanishes with mesh refinement, the solutions of all three meshes are nearly identical and agree well with the results in the literature.\cite{hauke1998comparative,xu2017compressible}

%\begin{figure}
%\begin{tikzpicture}
%\centering
%\begin{axis}[
%   	height=0.2\textheight, width=0.7\textwidth,
%	xmin = 0,
%    xmax = 1,
%    ymin = 0,
%    ymax = 1,
%    xlabel={Distance from leading edge, x},
%   	ylabel= {Pressure coefficient, $C_p$},
%   	grid=major,
%   	major grid style={dotted},
%   	   	legend style={
%		at={(0.5,0.8)},
%		anchor=north,
%		cells={anchor=west},
%		legend columns=0,
%	},
%   	legend entries={ Coarse, Medium, Fine},
%]    
%	\addplot+[no markers]  table [x index = 6, y index = 0,each nth point=10,col sep=comma] {../../CommonPics/fpCoarseNew.csv};
%	\addplot+[no markers,dashed]  table [x index = 6, y index = 0,each nth point=10,col sep=comma] {../../CommonPics/fpMediumNew.csv};
%	\addplot+[no markers]  table [x index = 6, y index = 0,each nth point=10,col sep=comma] {../../CommonPics/fpFineNew.csv};   
%\end{axis}
%\end{tikzpicture}
%\caption{Supersonic flow over a flat plate. Pressure coefficient along the solid wall.}
%\label{fig:fpPressure}
%\end{figure}

\subsection{Pressure pulse}
\label{ssec:pulse}
To evaluate the temporal accuracy of the three discretizations methods (FST, SST, UST), we derive the following pressure pulse test case in one space dimension. The test case takes into account the nonlinear characteristics of a sound wave\cite{thompson1971compressible} traveling in air at standard temperature and pressure (STP)
\begin{equation}
 p_0 = 10^5\, Pa, \quad T_0 = 273.15\, K, \quad c_0= \sqrt{\gamma R T_0} \approx 331.29 \, \frac{m}{s}.
\end{equation}

\begin{figure}
\centering
\includegraphics[width=0.9\textwidth]{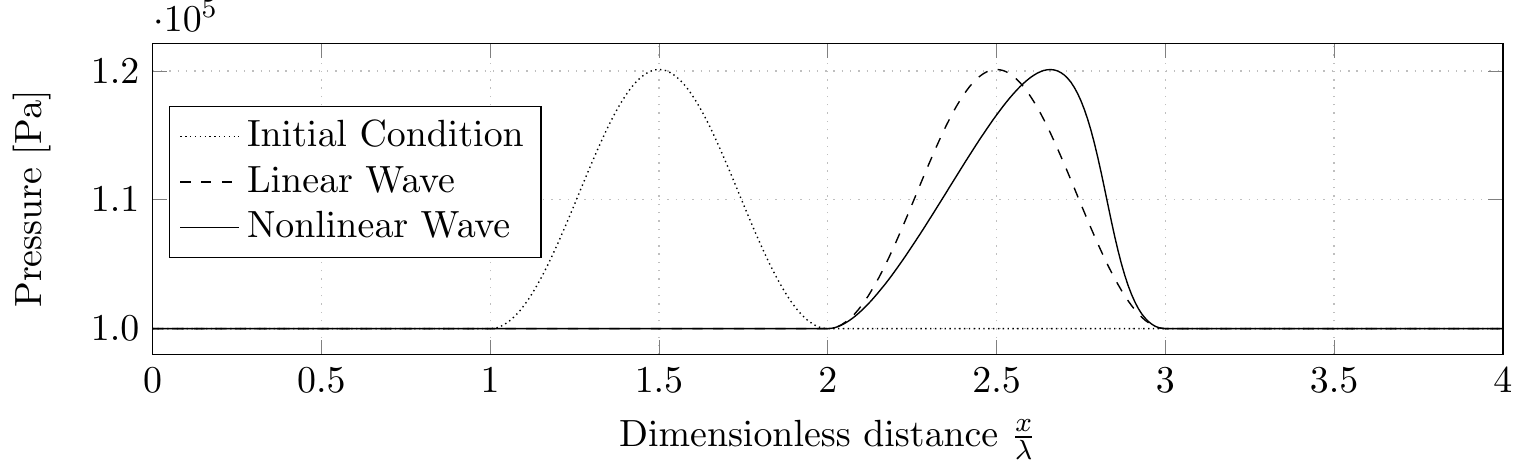}
\caption{Pressure pulse. Initial condition and waves after traveling one wave length at $t = t_f$.}
\label{fig:presPulse}
\end{figure}

As shown in Figure~\ref{fig:presPulse} the computational domain spans four wave lengths in $x$-direction, $0 \leq x \leq 4 \lambda$. By means of the Heaviside function, $\mathcal{H}(x)$, the perturbation function
\begin{equation}
S(x) = \alpha \left( 1- \cos\left(2\pi \frac{x}{\lambda}\right)\right) \cdot  \left( \mathcal{H}\left(\frac{x}{\lambda}-1\right)- \mathcal{H}\left(\frac{x}{\lambda}-2\right)\right),
\label{eq:perturb}
\end{equation}
defines the initial distribution of the speed of sound  as
\begin{equation}
c = c_0(1 + S(x)), \quad  \text{at} \quad t = 0.
\end{equation}
Using isentropic relations, initial conditions for the pressure-primitive variables are defined as
\begin{eqnarray}
p &= p_0 \left( \frac{c}{c_0}\right)^{\frac{2\gamma}{\gamma-1}} &= \left(1+S(x)\right)^{\frac{2\gamma}{\gamma-1}} p_0,  \\
u_1 &= \frac{2}{\gamma-1}(c-c_0) &= \frac{2}{\gamma-1} S(x) c_0,\\
u_2 & &= 0,\\
T & = T_0 \left( \frac{c}{c_0}\right)^2 & = \left(1+S(x)\right)^2 T_0.
\end{eqnarray}

After the time period $t_f = \frac{\lambda}{c_0}$, the head and tail of the wave have traveled exactly one wave length $\lambda$. Using the theory on shock formation in one-dimensional unsteady flows,\cite{thompson1971compressible} we compute the time $t^*$ after which a shock first forms as
\begin{equation}
t^* = - \frac{\gamma-1}{\gamma+1} \left( \min\left(\ddx{c}{x}\right) \right)^{-1} = \frac{\gamma-1}{\gamma+1} \frac{1}{2\pi\alpha} \frac{\lambda}{c_0}.
\end{equation}
Therein, $\min\left(\ddx{c}{x}\right)$ corresponds to the steepest descent along the wave. For our test case, we further request $ t_f = \frac{1}{2} t^*$ to ensure that the nonlinear behavior is clearly visible, but the gradients of the solution remain finite. This condition yields $\alpha = \frac{1}{24\pi}$ in the perturbation function $S(x)$ in Equation~\eqref{eq:perturb}.

At $x = 0$ and $x = 4 \lambda$,  the conditions of the fluid at rest are prescribed as Dirichlet boundary conditions
\begin{equation}
	     \bY_{DBC} = \left[  \begin{array}{c}
	     p \\
	     u_1\\
	     u_2\\
	     T \\
	      \end{array} \right] = \left[  \begin{array}{c}
	     p_0 \\
	     0 \\
	     0 \\
	     T_0 \\
	      \end{array} \right].
\end{equation}

For an inviscid fluid, the nonlinear acoustics theory can be used to obtain a reference solution.\cite{thompson1971compressible} The reference solution is constructed by advancing the all points of the wave by $(u+c) t_f$. Note that $u$ and $c$ vary along $x$. Figure~\ref{fig:presPulse} displays the initial condition for the pressure $p(x, t=0)$ and compares the nonlinear reference solution $p(x,t = t_f)$ to the linear case.

\begin{figure}
\centering
\includegraphics[width=0.9\textwidth]{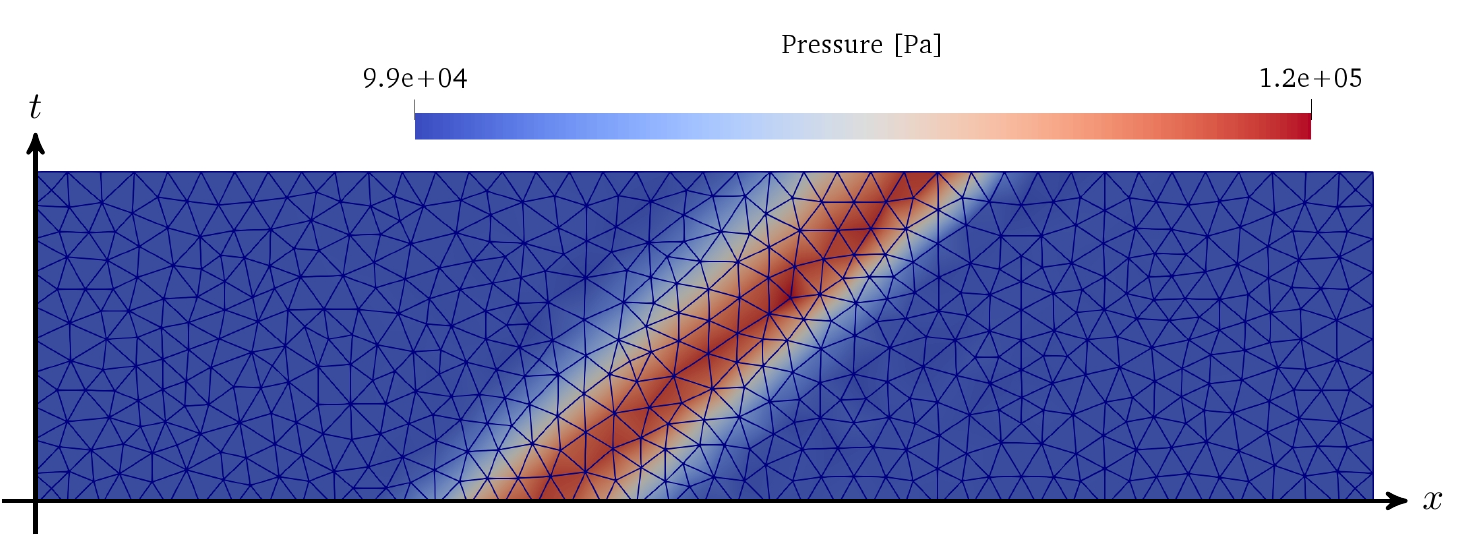}
\caption{Pressure pulse. Solution on coarse UST mesh.}
\label{fig:presUST}
\end{figure}

To give a first example of a solution on an unstructured space-time (UST) mesh, Figure~\ref{fig:presUST} shows the solution obtained on an extremely coarse UST mesh. The space dimension is on the horizontal axis and the time dimension on the vertical axis.

Lacking an implementation of the space-time methods for a single space dimension, we discretize the quasi one-dimensional spatial computational domain ($0 \leq x \leq 4 \lambda$, $0 \leq y \leq \frac{\lambda}{100}$) by splitting $400 \times 1$ square elements into right triangles. At $y=0$ and $y=\frac{\lambda}{100}$, the discretization has 100 elements per wave length. For all three space-time methods, the spatial discretization of the initial and final time level are identical. The temporal discretization size $\Delta t$ is chosen such that a CFL-number of $C_{\Delta t} = \frac{\Delta t c_0}{\Delta x} = 2$ is obtained and the temporal discretization error dominates the spatial discretization error.

\begin{figure}
\centering
\includegraphics[width=0.9\textwidth]{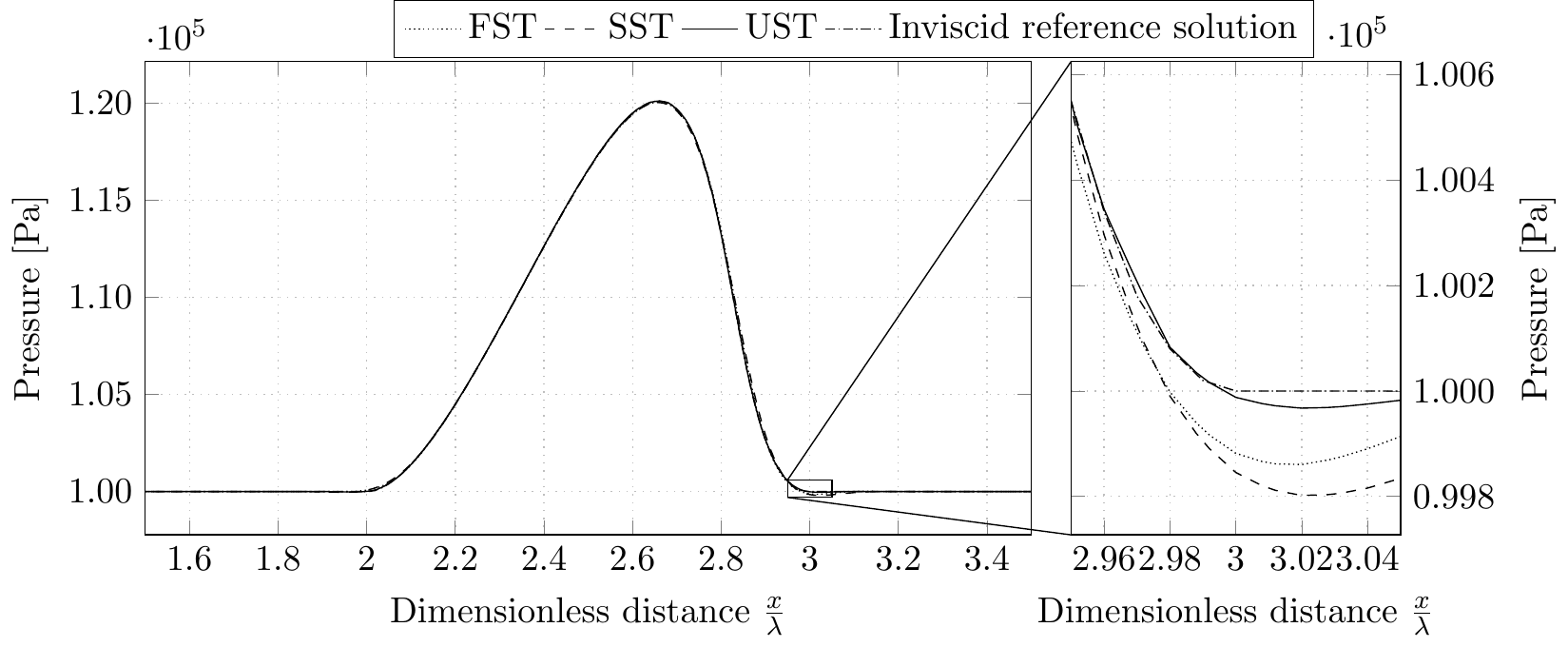}
\caption{Pressure pulse. Solutions at $t = t_f$ for a CFL-number of $C_{\Delta t} = \frac{\Delta t c_0}{\Delta x} = 2$}
\label{fig:presSol}
\end{figure}

Looking at the complete wave after travelling one wave length (Figure~\ref{fig:presSol} on the left), the solutions on all three space-time discretizations coincide with the inviscid nonlinear reference solution. Zooming in on the head of the wave (Figure~\ref{fig:presSol} on the right), a small undershoot of approximately 0.2\% is visible for the FST and SST solution. The undershoot of the UST solution is significantly smaller (about 0.04\%). These results illustrate the good temporal accuracy of all three space-time discretization methods.

\subsection{Valve}
\label{ssec:valve}
\begin{figure}
\centering
\includegraphics[width=\textwidth]{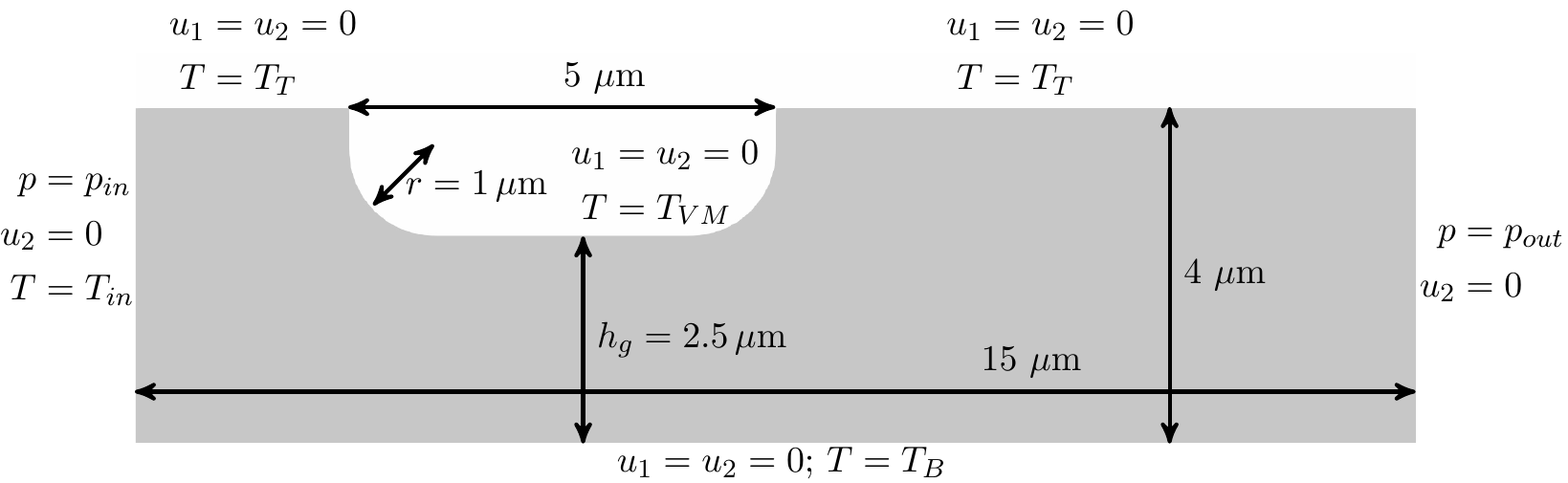}
\caption{Valve. Fluid domain in open configuration with boundary conditions}
\label{fig:valveSetup}
\end{figure}

The UST method allows to solve the compressible Navier--Stokes equations on spatial computational domains that undergo topological changes over time. Namely, with space-time finite elements one can discretize the space-time continuum of a spatial computational domain, which splits into two and reconnects over time, conformingly. This is demonstrated here using a two-dimensional valve test case in which the fluid domains are separated by a valve member. The dimensions and conditions of this test case are inspired by the gas flow through the narrow gaps in the piston ring pack of an internal combustion engine (see Section~\ref{ssec:pistonRings}). The problem setup is displayed in Figure~\ref{fig:valveSetup}. The fluid domain of the valve consist of a two-dimensional channel of 15 $\mu$m length and 4 $\mu$m height and a valve member that creates a 2.5 $\mu$m gap with the valve bottom in the open configuration. In the closed configuration, the gap between the valve member and the channel bottom vanishes completely ($h_g=0$).

\begin{figure}
\centering
\includegraphics[width=0.4\textwidth]{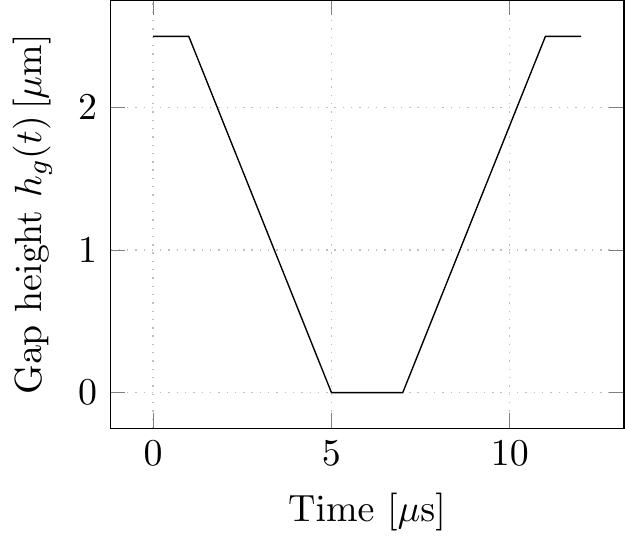}
\caption{Valve. Gap height during valve cycle}
\label{fig:valveHg}
\end{figure}

\begin{figure}
\centering
\includegraphics[width=0.5\textwidth,trim={0cm 0cm 0cm 0cm},clip]{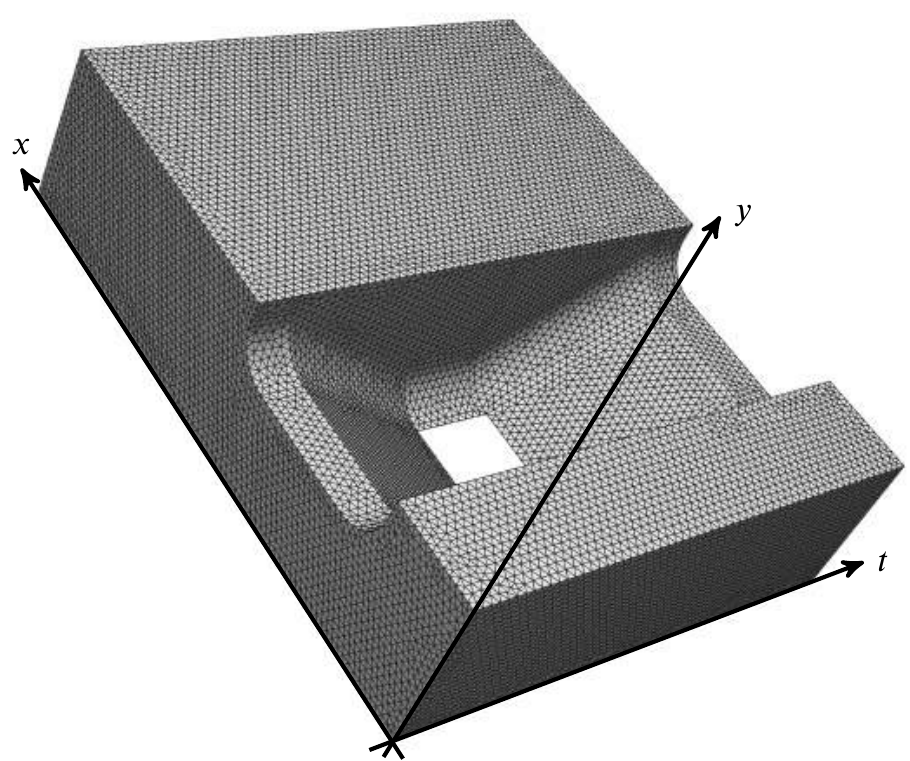}
\caption{Valve. Unstructured space-time mesh.}
\label{fig:valveMesh}
\end{figure}

The boundary conditions are as follows. On the left domain boundary, the pressure $p_{in}= 1.2\times10^5 Pa$ and the temperature $T_{in} = 373.15 K$ is prescribed; the tangential component of the velocity is set to zero $u_2 = 0 \frac{m}{s}$. On the right domain boundary, the pressure $p_{out}= 1\times10^5 Pa$ is prescribed and the tangential component of the velocity is set to zero $u_2 = 0 \frac{m}{s}$. On these two open boundaries, the viscous contributions ($\tau, q \neq 0$) are included in the boundary integral $\int_{P_h} \bW \cdot \bh \, dP$. On the solid walls of the channel bottom, top and the valve member, no-slip boundary conditions are applied and the temperatures $T_{B} = 393.15 K$, $T_{T} = 383.15 K$, and $T_{VM} = 403.15 K$ are prescribed, respectively.

The temporal evolution of the valve cycle is characterized by the gap height $h_g(t)$ between valve member and channel bottom plotted in Figure~\ref{fig:valveHg}. Initially, the valve is held open for 1 $\mu$s, then closed over the following 4 $\mu$s. Subsequently, it is completely shut for 2 $\mu$s, before it is again opened over the next 4 $\mu$s, and finally held open for another 1 $\mu$s. The resulting space-time domain is discretized with the UST method, yielding the mesh shown in Figure~\ref{fig:valveMesh}.

In the open valve configuration, the pressure gradient of 0.2 bar leads to a laminar flow with a maximum Reynolds number $Re = \frac{\rho\, u_{max}\, h_g(12\, \mu s)}{\mu_{ref}} \approx 10$ and Mach number $Ma = \frac{u_{max}}{\sqrt{\gamma R T}} \approx 0.21$. The fluid density varies between $0.88 \frac{kg}{m^3}$ and $1.12 \frac{kg}{m^3}$. Figure~\ref{fig:valveT} shows the pressure, velocity, and temperature distribution in the closed, half open, and fully opened valve. The renderings are obtained by evaluating the space-time data on the $x$-$y$-planes at $t=6, \, 9, \, 12\, \mu$s. In the closed configuration, the valve member completely separates the fluids to its left and right. The gas at the left attains the high pressure value prescribed at the inlet, the gas at the right the low pressure value prescribed at the outlet (Figure~\ref{fig:valveP6}). As the fluid is at rest (Figure~\ref{fig:valveV6}), the temperature distribution is governed by the diffusion between the prescribed values at the inlet and the solid walls (Figure~\ref{fig:valveT6}). In the open configurations, the pressure expands from the inlet on the left through the valve to the outlet on the right  (Figure~\ref{fig:valveP9} and~\ref{fig:valveP12}). The fluid velocities for the fully opened configuration (Figure~\ref{fig:valveV12}) are significantly higher than the fluid velocities in the half opened configuration (Figure~\ref{fig:valveV9}). Figure~\ref{fig:valveT12} clearly shows the influence of the fluid flow on the temperature distribution in the valve. In summary, our valve simulation shows, that the UST method is capable to simulate flows on spatial domains with time-varying topology.
 
\newcommand{\myW}{0.33\textwidth} 
\begin{figure}
\centering
\includegraphics[width=0.32\textwidth,trim={16cm 0cm 18cm 14cm},clip]{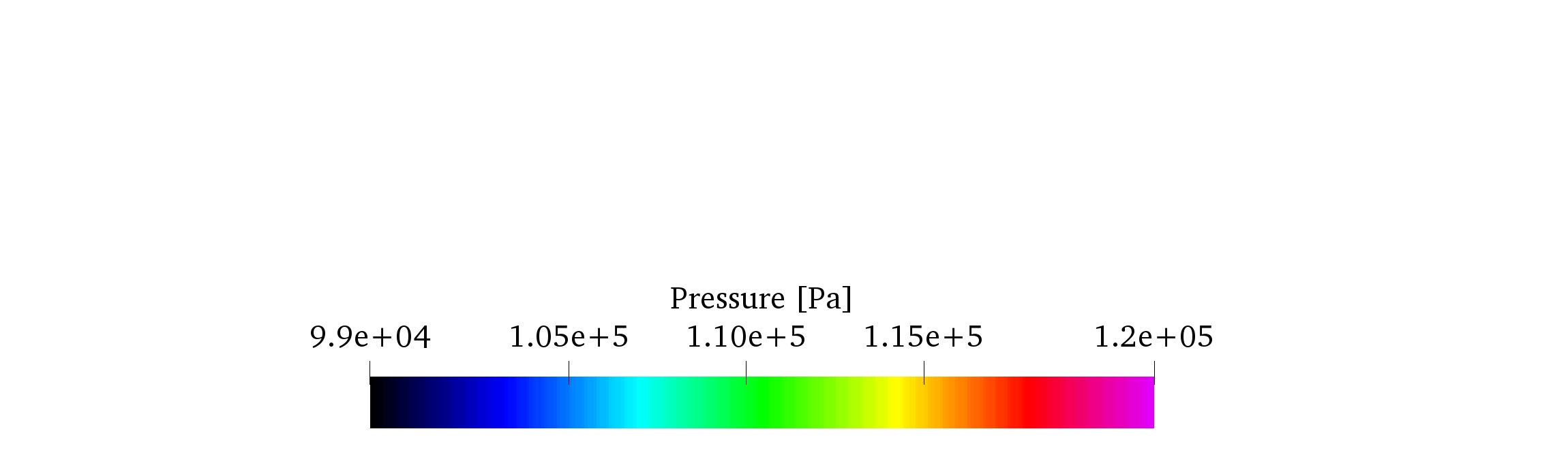}
\includegraphics[width=0.32\textwidth,trim={16cm 0cm 18cm 14cm},clip]{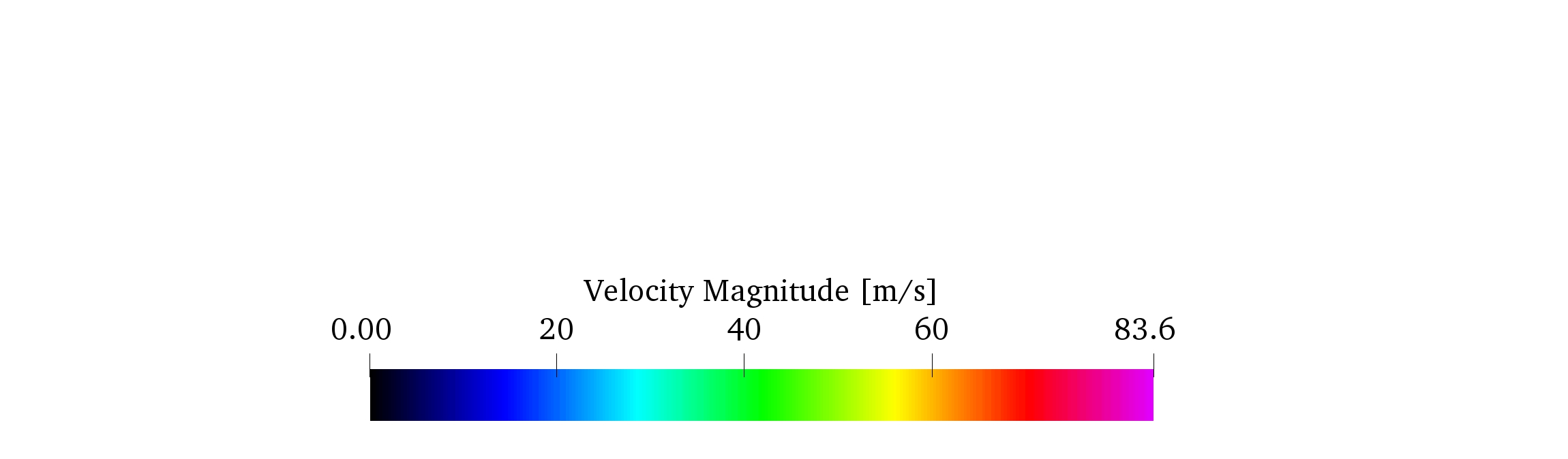}
\includegraphics[width=0.32\textwidth,trim={16cm 0cm 18cm 14cm},clip]{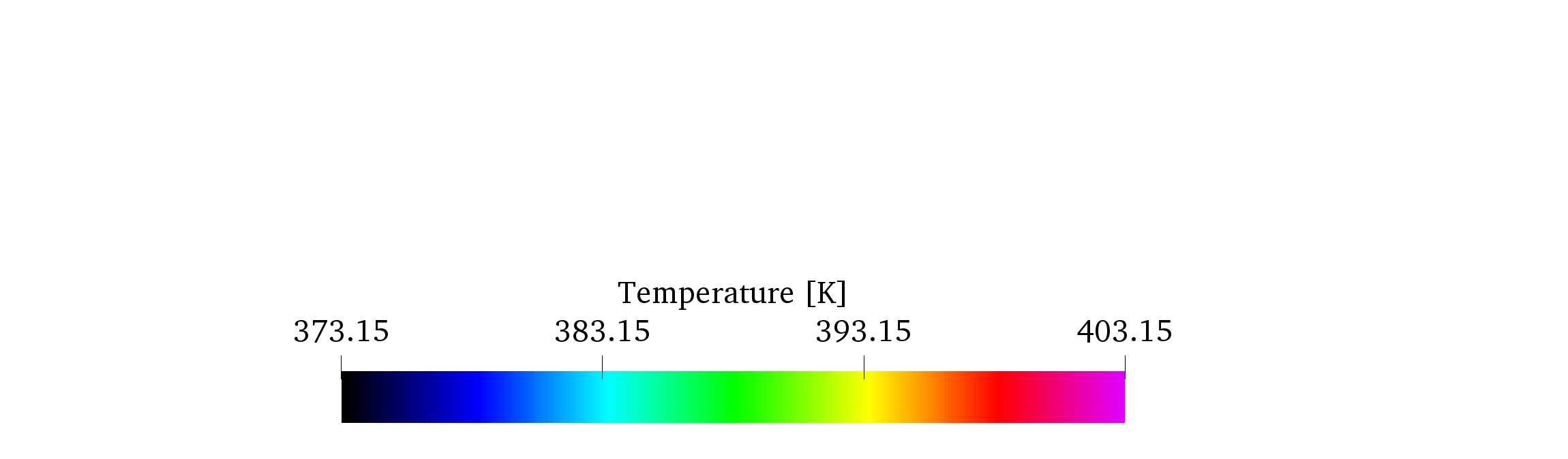}

\subfloat[Pressure distribution at $t=6\, \mu$s \label{fig:valveP6}]{\includegraphics[width=\myW,trim={0cm 0cm 0cm 0cm},clip]{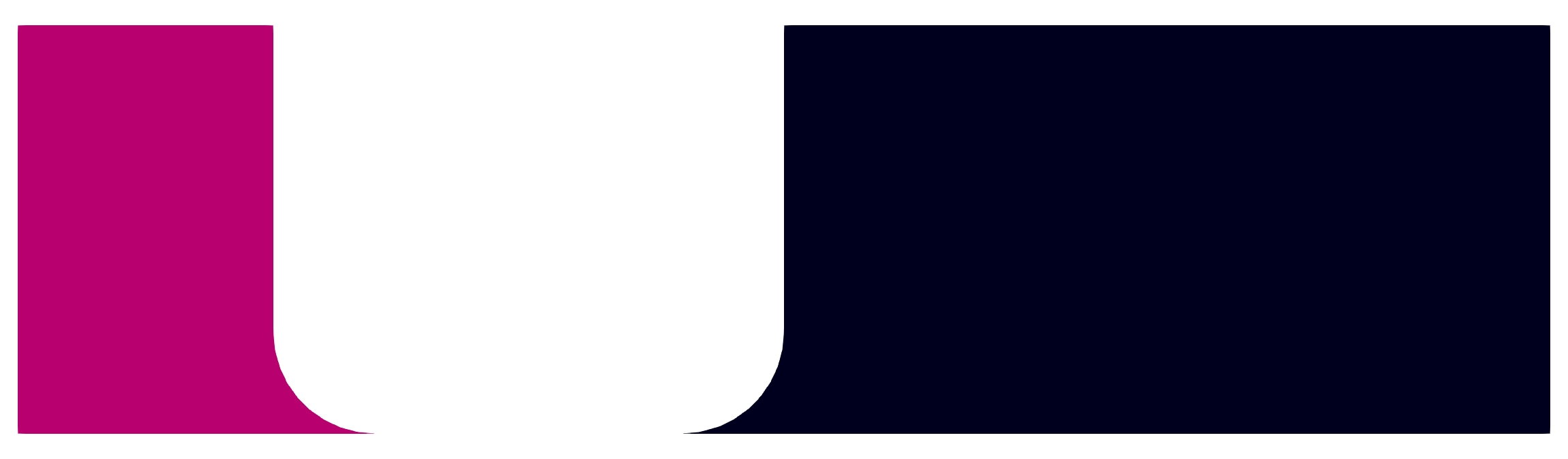}}
\subfloat[Velocity field at $t=6\, \mu$s \label{fig:valveV6}]{\includegraphics[width=\myW,trim={0cm 0cm 0cm 0cm},clip]{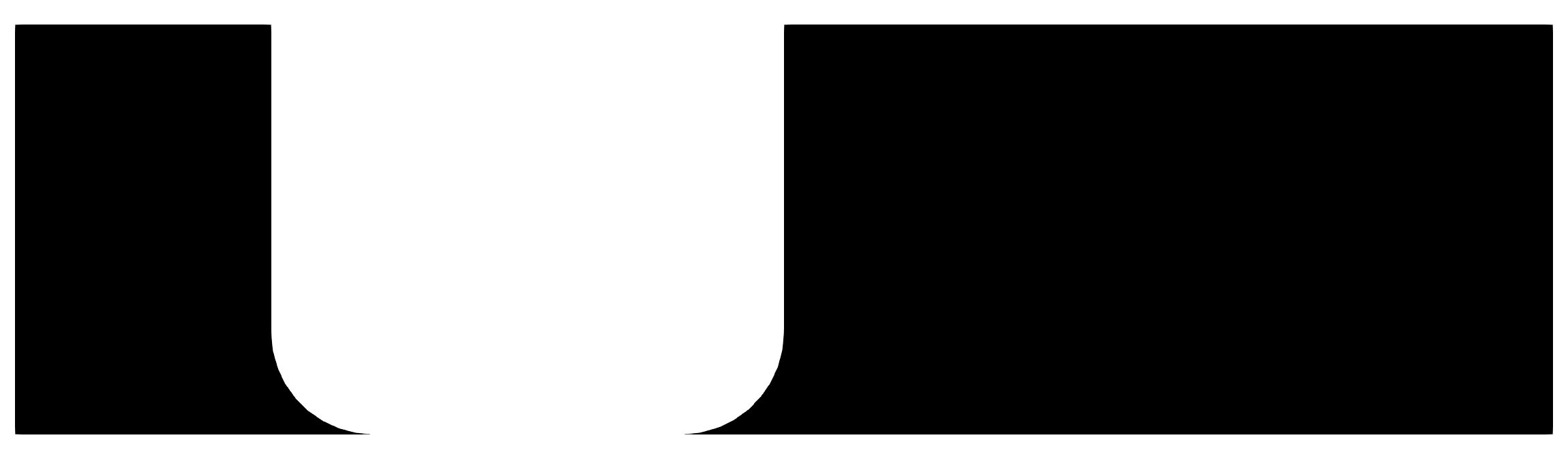}}
\subfloat[Temperature distribution at $t=6\, \mu$s \label{fig:valveT6}]{\includegraphics[width=\myW,trim={0cm 0cm 0cm 0cm},clip]{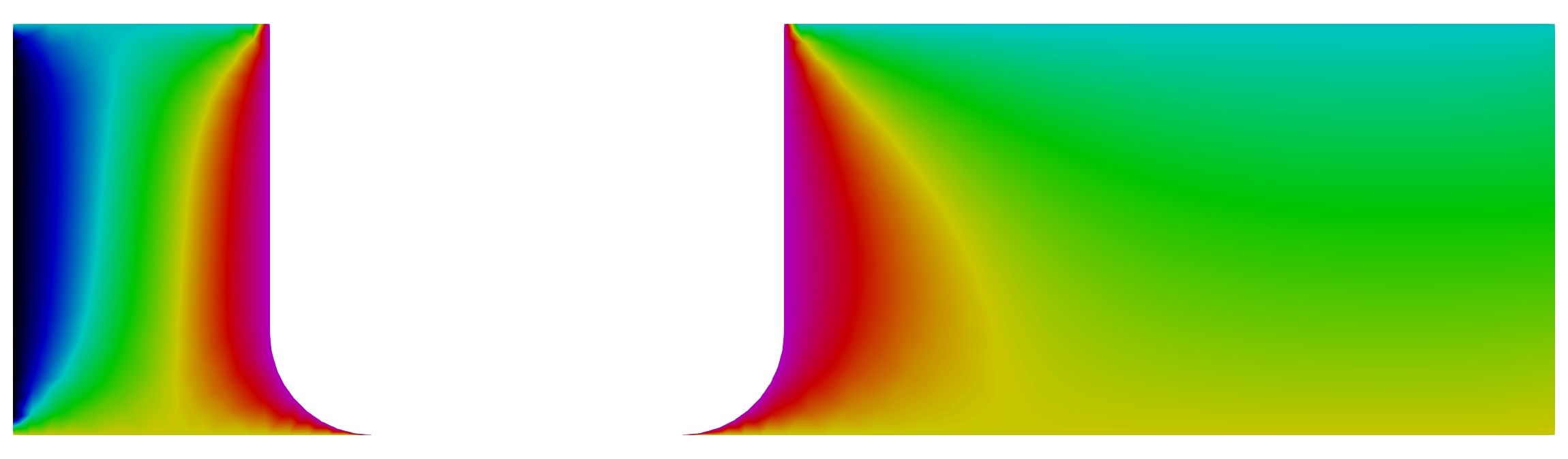}}

\subfloat[Pressure distribution at $t=9\, \mu$s \label{fig:valveP9}]{\includegraphics[width=\myW,trim={0cm 0cm 0cm 0cm},clip]{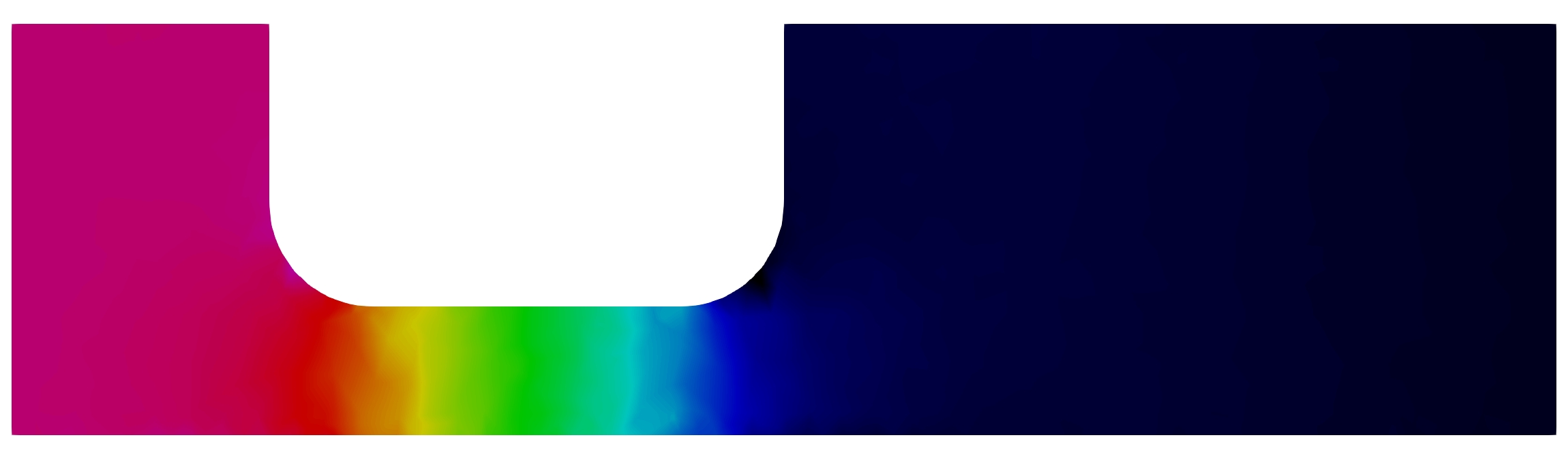}}
\subfloat[Velocity field at $t=9\, \mu$s \label{fig:valveV9}]{\includegraphics[width=\myW,trim={0cm 0cm 0cm 0cm},clip]{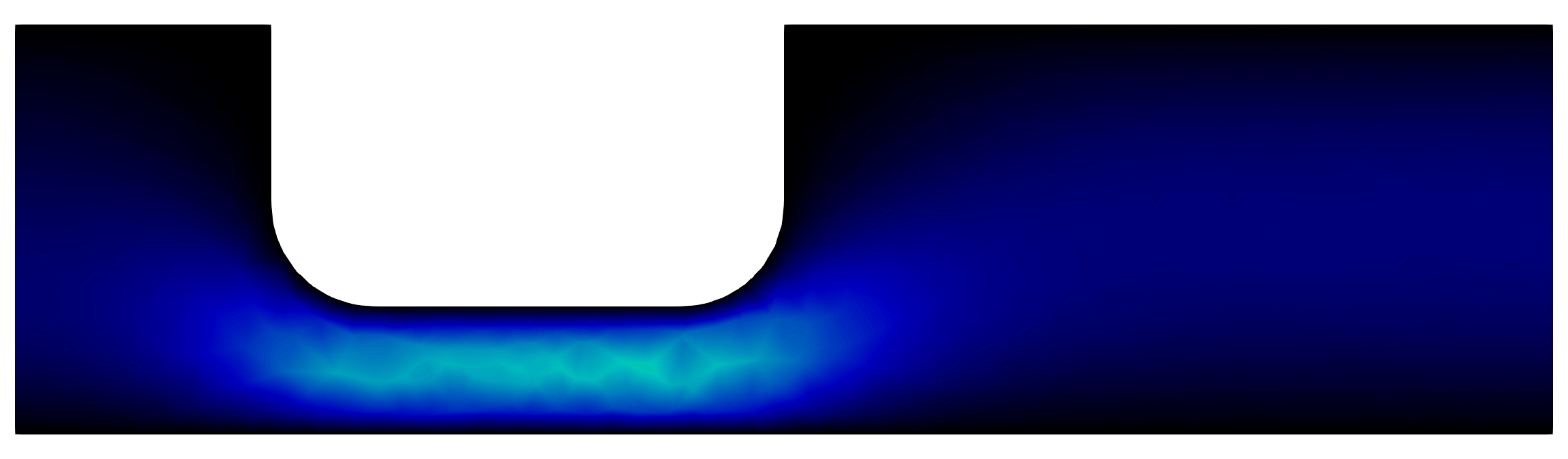}}
\subfloat[Temperature distribution at $t=9\, \mu$s \label{fig:valveT9}]{\includegraphics[width=\myW,trim={0cm 0cm 0cm 0cm},clip]{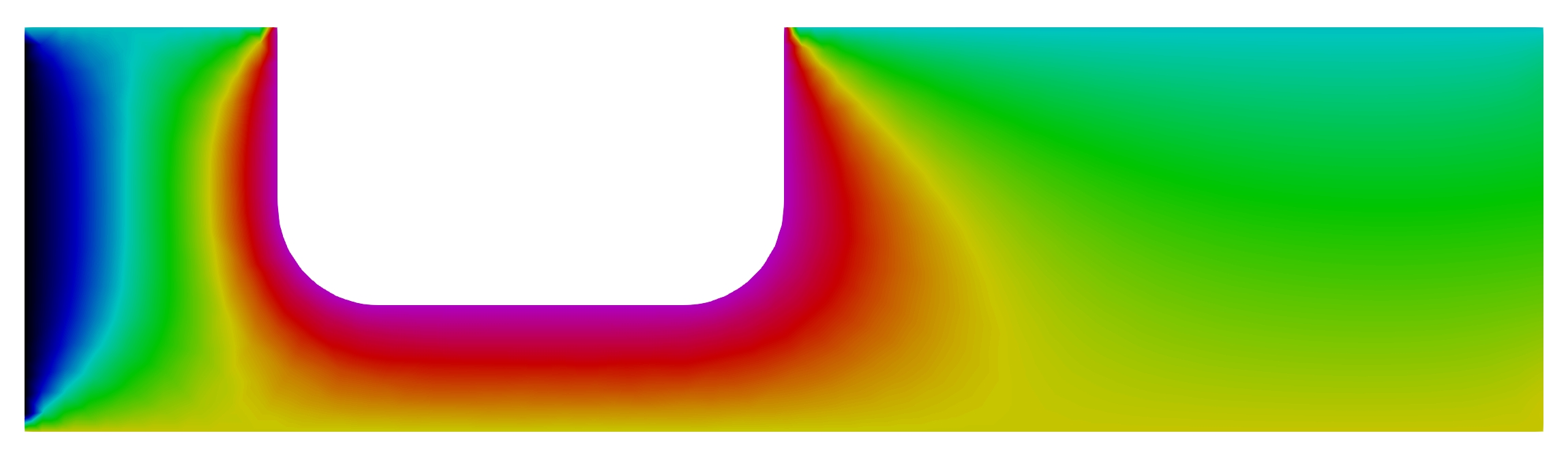}}\\

\subfloat[Pressure distribution at $t=12\, \mu$s \label{fig:valveP12}]{\includegraphics[width=\myW,trim={0cm 0cm 0cm 0cm},clip]{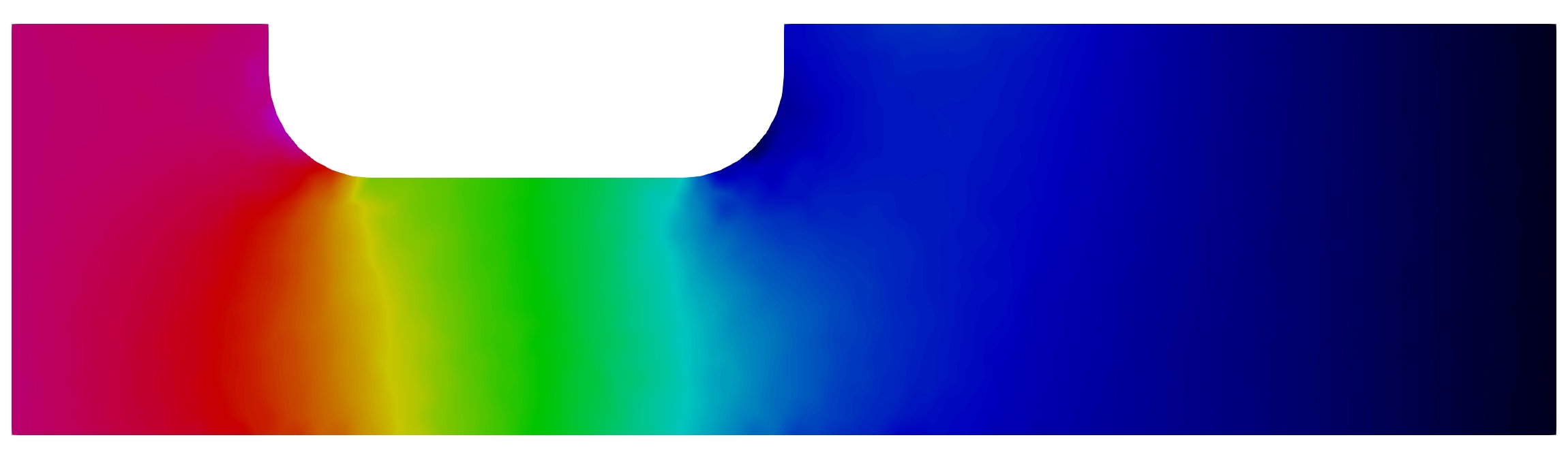}}
\subfloat[Velocity field at $t=12\, \mu$s \label{fig:valveV12}]{\includegraphics[width=\myW,trim={0cm 0cm 0cm 0cm},clip]{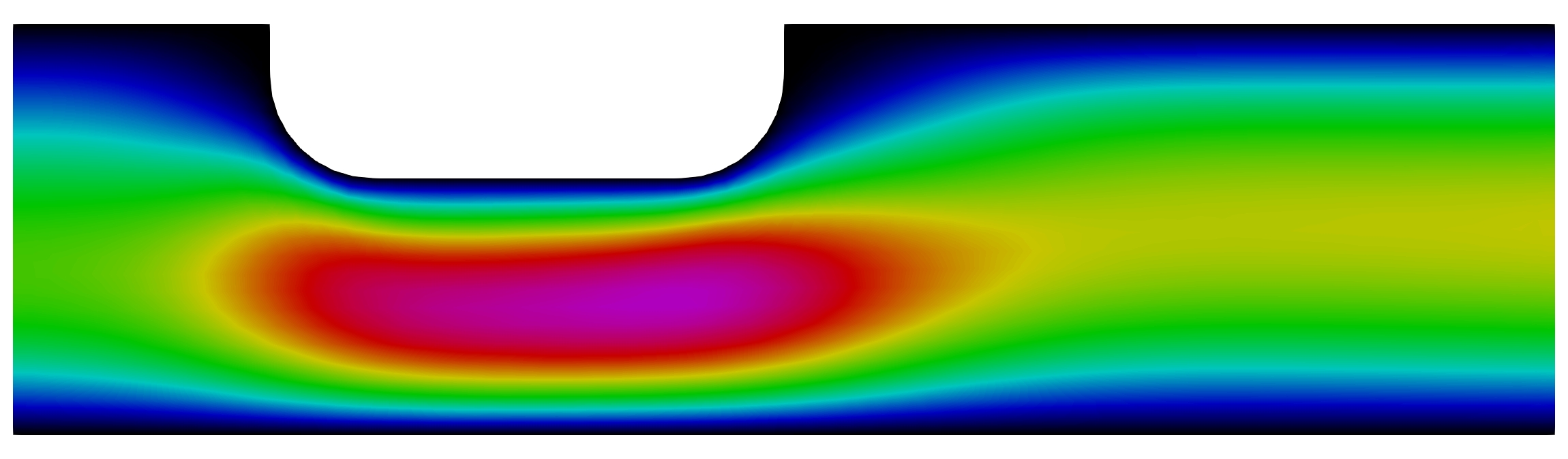}}
\subfloat[Temperature distr. at $t=12\, \mu$s \label{fig:valveT12}]{\includegraphics[width=\myW,trim={0cm 0cm 0cm 0cm},clip]{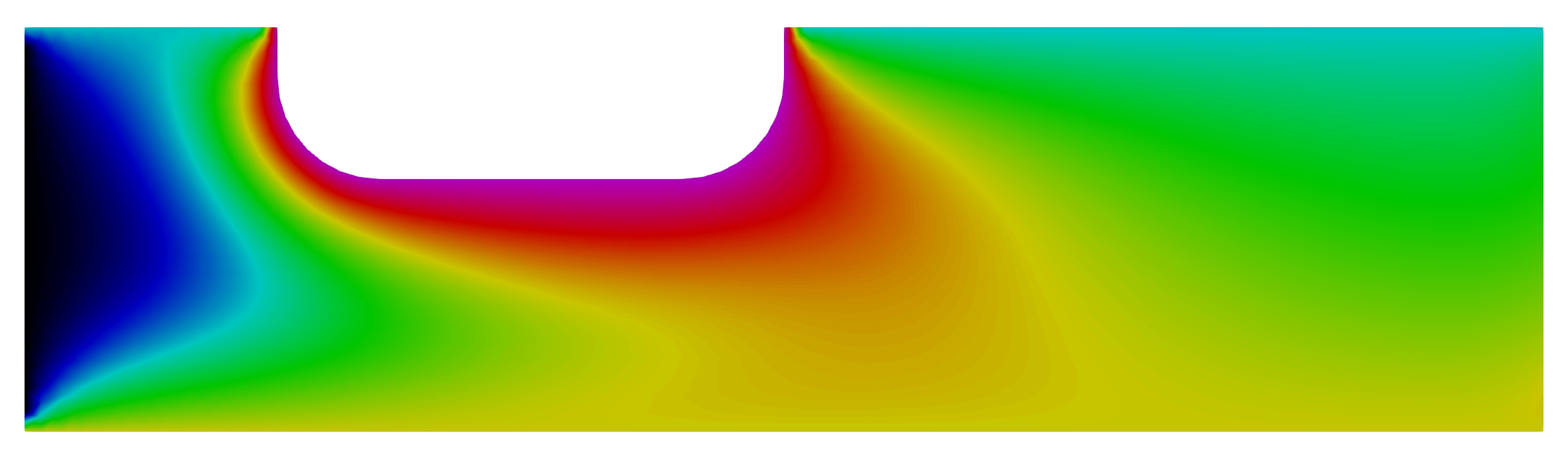}}\\
\caption{Pressure, velocity, and temperature distribution in the closed, half open, and fully opened valve.}
\label{fig:valveT}
\end{figure} 

\subsection{Blow-by past piston rings}
\label{ssec:pistonRings}
This test case derives from a problem of engineering interest and demonstrates the potential to save computation time by local temporal refinement in simplex space-time meshes. As the fuel efficiency of an internal combustion engine is directly related to the performance of the
piston ring pack, simulating the blow-by past piston rings is an active field of research.\cite{oliva2016numerical} Figure~\ref{fig:pistonpic} relates the considered computational domain to the piston ring pack of an exemplaric piston. The domain for the gas flow computation consists of the voids in the piston ring pack. The geometry of the computational domain is simplified, such that the three piston lands are only connected through the ring end gaps. Still, scales that need to be resolved range from $\approx0.1$ mm in the ring end gaps to $\approx10$ cm (piston diameter) making the mesh generation a intricate task. The tetrahedral discretization of the spatial computaitonal domain (Figure~\ref{fig:pistonpic}) is strongly adapted to the expected flow field, which leads to highly stretched elements. The mesh was generated with GMSH.\cite{geuzaine2009gmsh}

\newcommand{\myTW}{0.48\textwidth}
\begin{figure}
\centering
\begin{tabular}{cc}
\multirow[c]{2}{*}[2.2cm]{\subfloat[Computational domain\label{fig:pistonpic}]{\includegraphics[width=\myTW]{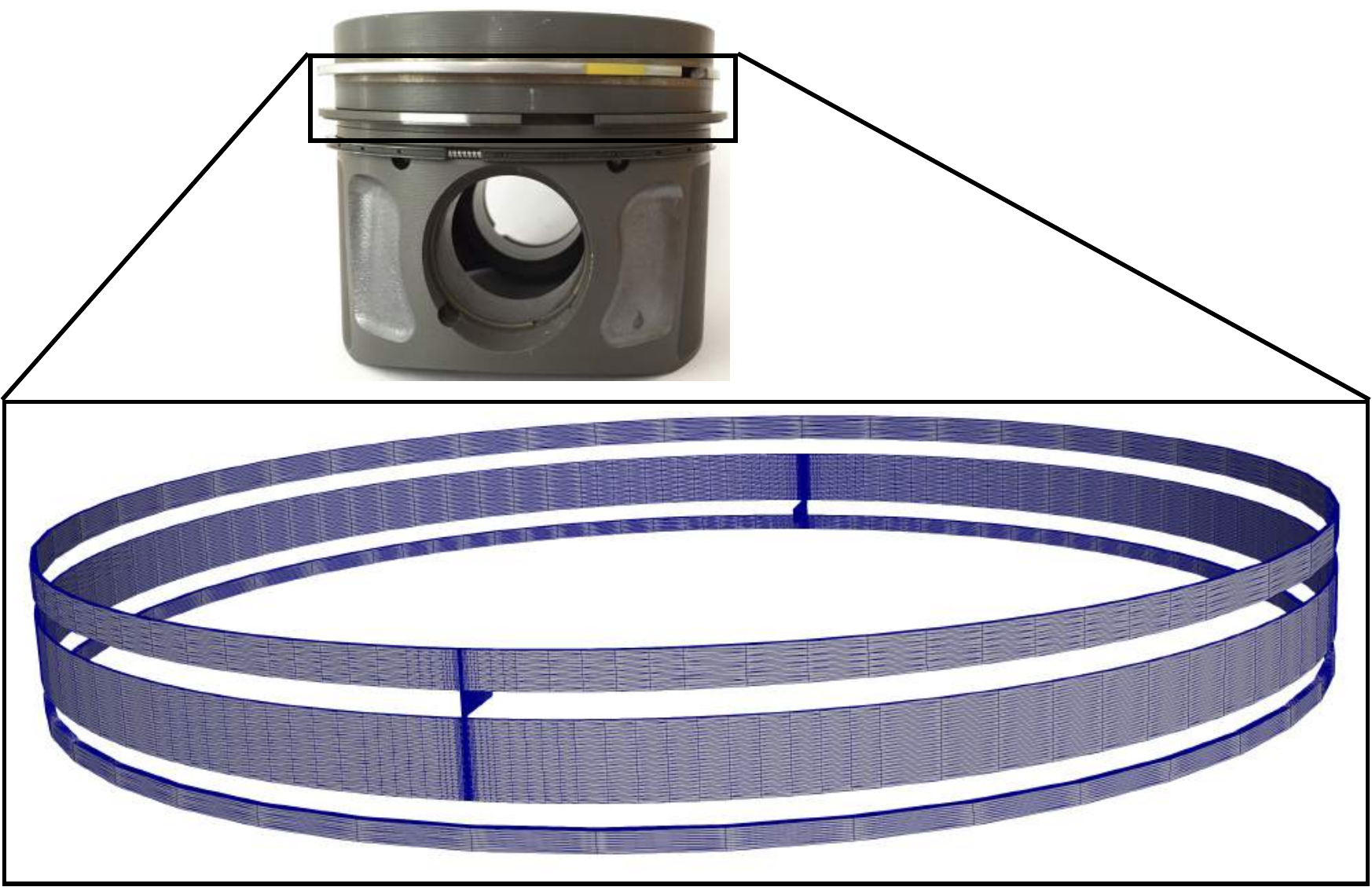}}}&
\subfloat[Pressure field after one degree crank angle\label{fig:RingPFull}]{\includegraphics[width=\myTW,trim={15cm 18cm 0cm 30cm},clip]{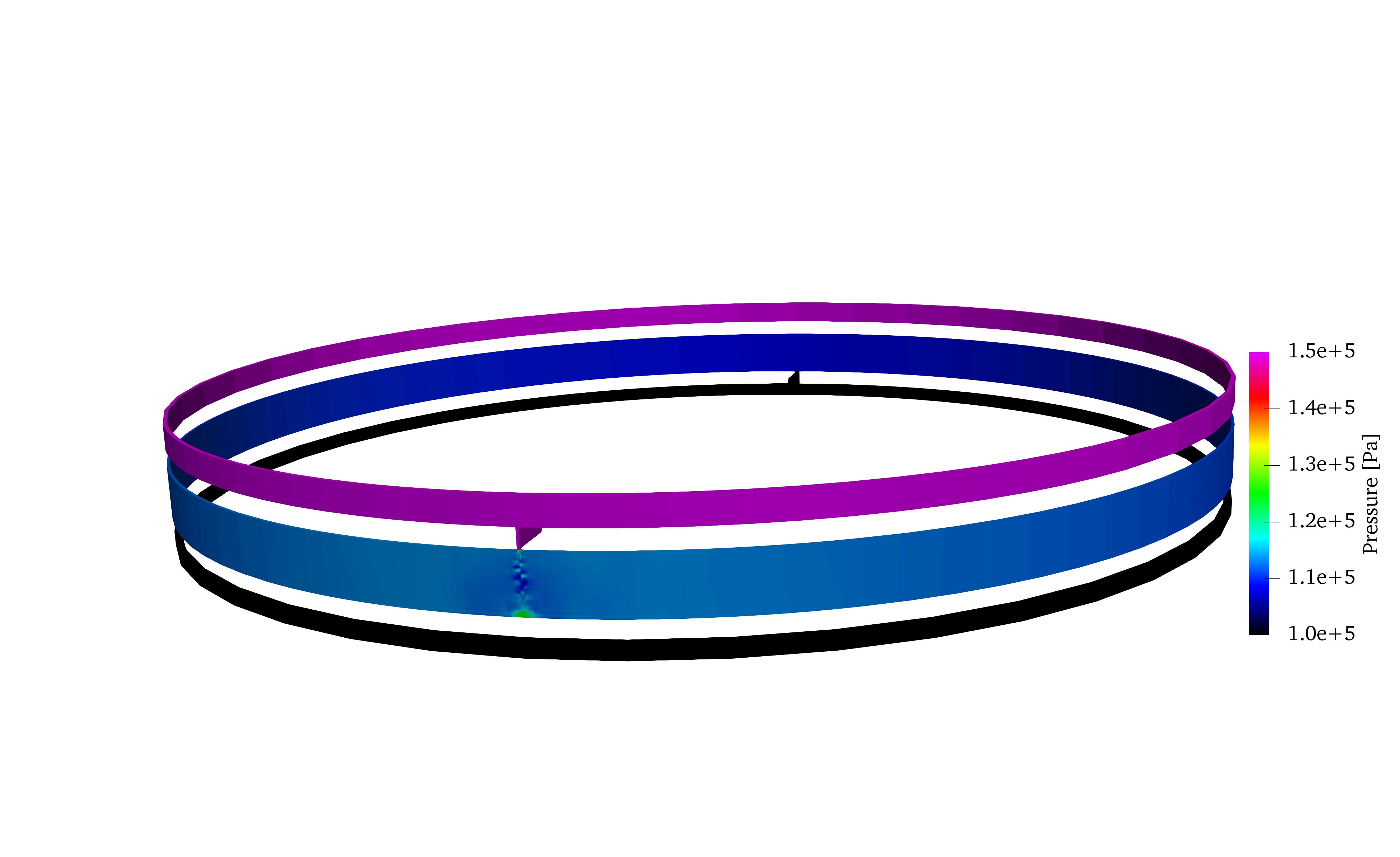}}\\
&
\subfloat[Velocity magnitude along streamlines after one degree crank angle\label{fig:RingVFull}]{\includegraphics[width=\myTW,trim={15cm 18cm 0cm 30cm},clip]{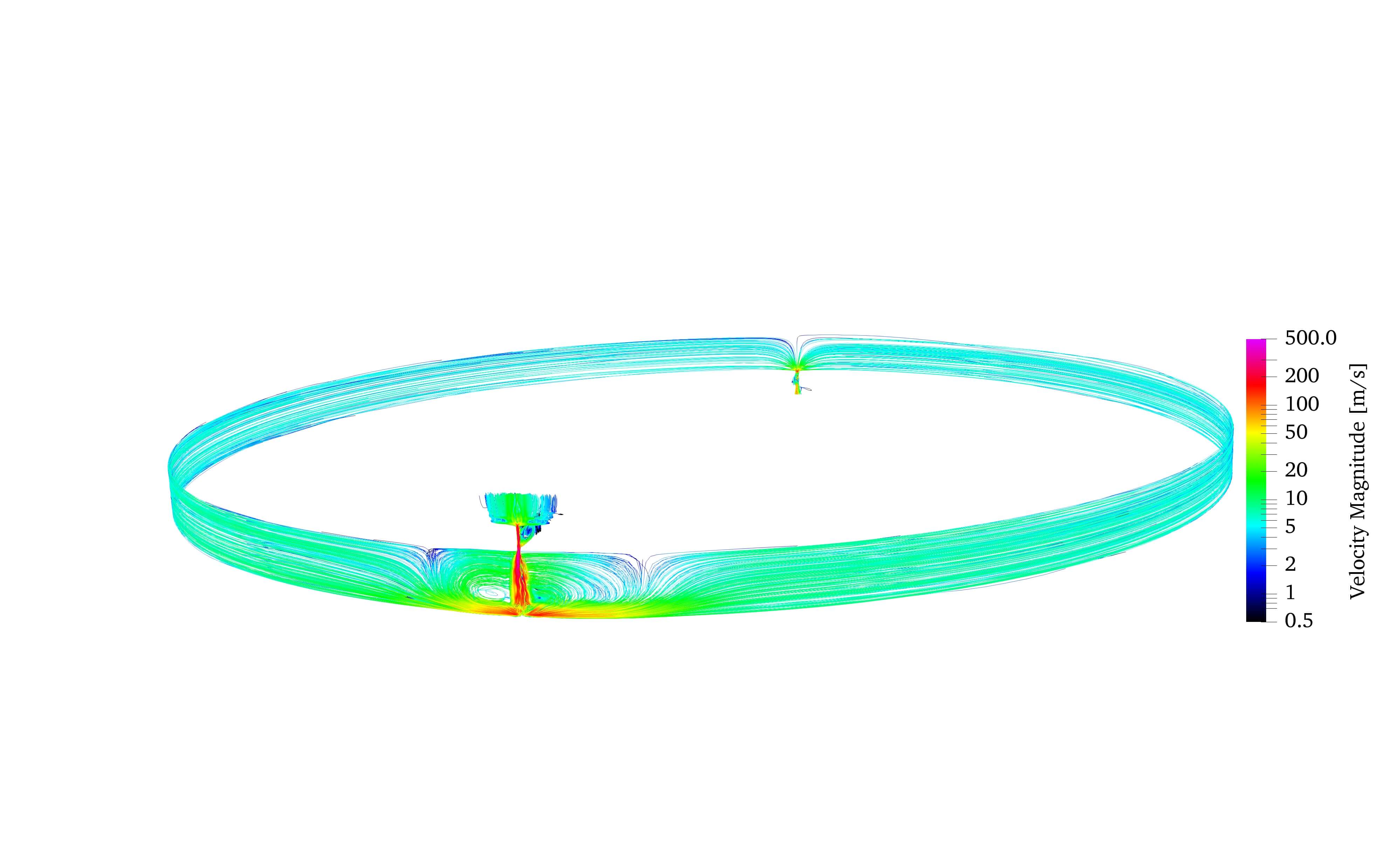}}\\
\subfloat[Area of four-dimensional temporal refinement\label{fig:RingTempRef}]{\includegraphics[width=\myTW,trim={15cm 18cm 0cm 30cm},clip]{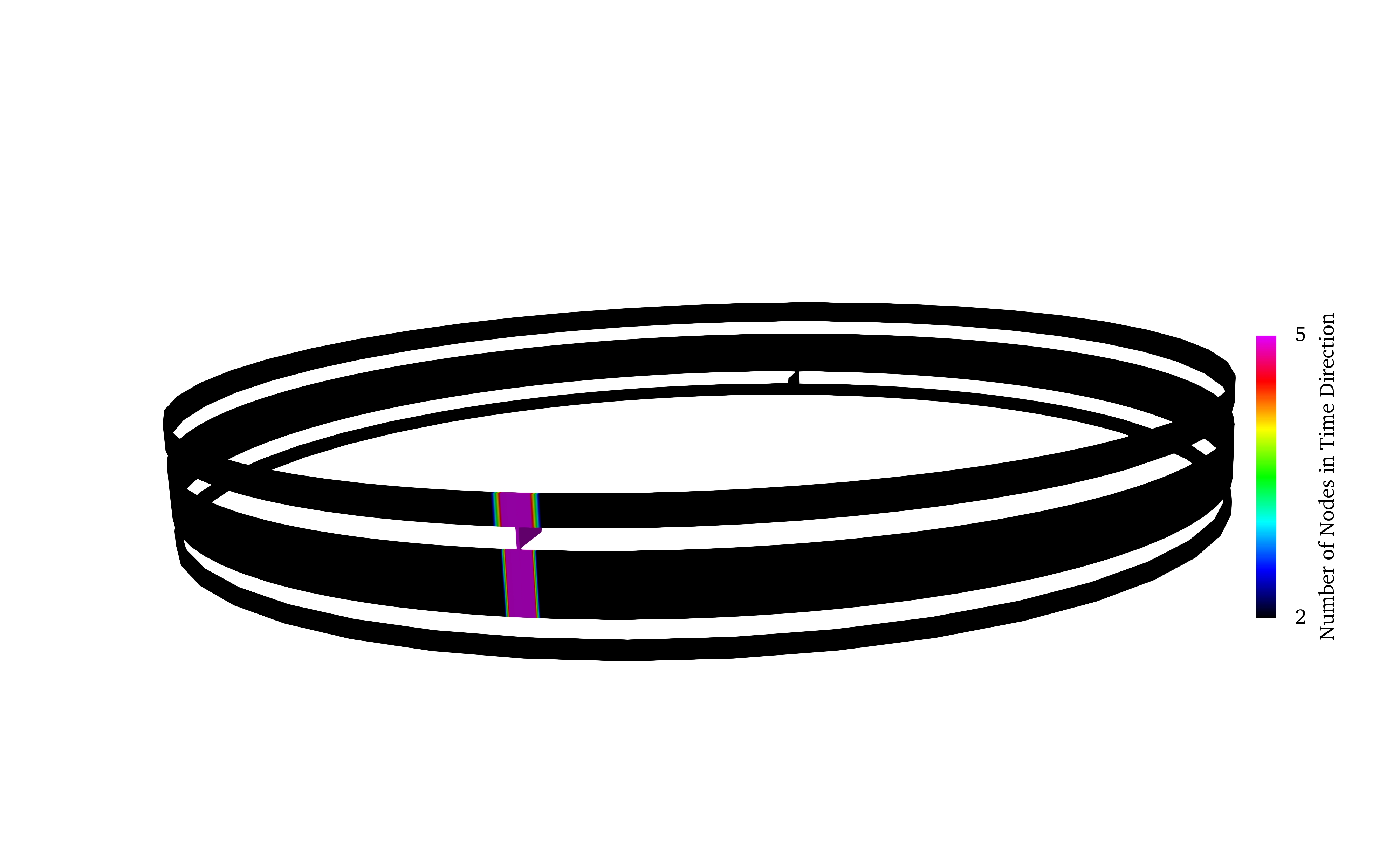}}&
\subfloat[Temperature along streamlines after one degree crank angle\label{fig:RingTFull}]{\includegraphics[width=\myTW,trim={15cm 18cm 0cm 30cm},clip]{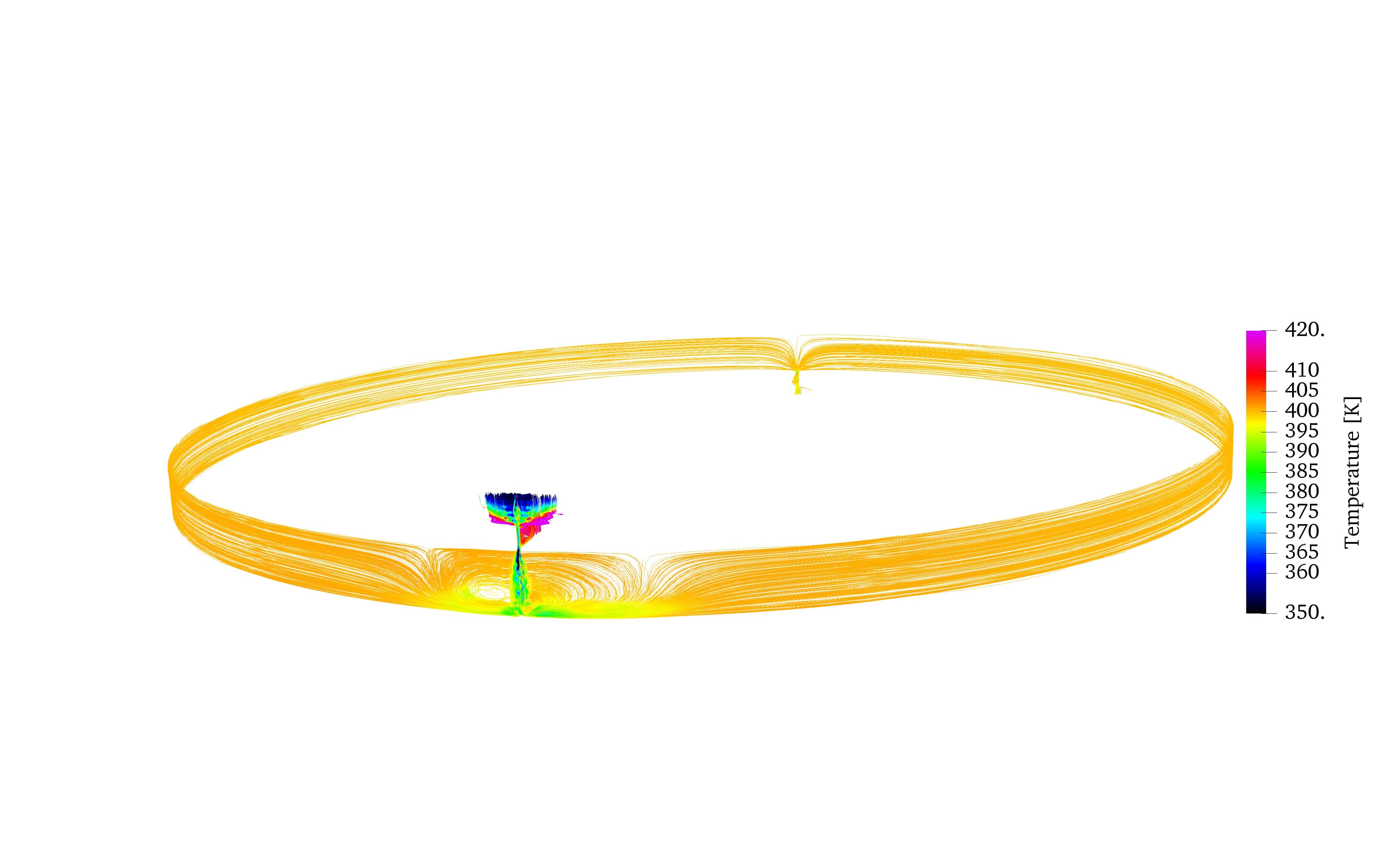}}
\end{tabular}
\caption{Simulation of blow-by in the piston ring pack.}
\label{fig:RingFull}
\end{figure}

The following boundary conditions characterize this test case. On the uppermost plane of the computational domain, the combustion chamber pressure during intake of $p_{in} = 2.174\times10^5 Pa$ and a temperature of $T_{in} = 352.3 K$ are applied. The crank case conditions of $p_{out} = 1\times10^5 Pa$ and $T_{out} = 400 K$ are prescribed on the lowest plane of the computational domain. On all solid walls no-slip conditions and a wall temperature of $T_W = 400 K$ are prescribed. A viable initial condition is obtained by gradually increasing the pressure at the inflow until the desired value of $p_{in}$ is reached. 

At an engine speed of 6000 rpm, one degree crank angle is simulated. We compare three different time discretizations. First, we consider an FST discretization which has 100 time steps per degree crank angle, opposed to an FST discretization with 400 time steps. Additionally, we apply the time refinement strategy of Behr\cite{behr2008simplex} to generate a four-dimensional SST discretization. The SST discretization has a temporal resolution that corresponds to 400 time steps per degree crank angle in the expected peak velocity region, while the resolution in the remaining spatial domain corresponds to 100 time steps per degree crank angle (see Figure~\ref{fig:RingTempRef}). Note that the spatial discretization for all three simulations is identical.

The flow field after one degree crank angle computed with 400 time steps on the FST discretization is shown in Figures~\ref{fig:RingPFull},\ref{fig:RingVFull}, and \ref{fig:RingTFull}. Vortical structures and peak velocities are restricted to the vicinity of the ring end gaps. Away from the ring end gaps, the gas flow along the second piston land is approximately uniform in circumferential direction.

\begin{figure}
\centering
\includegraphics[width=\textwidth,trim={3cm 3cm 3cm 70cm},clip]{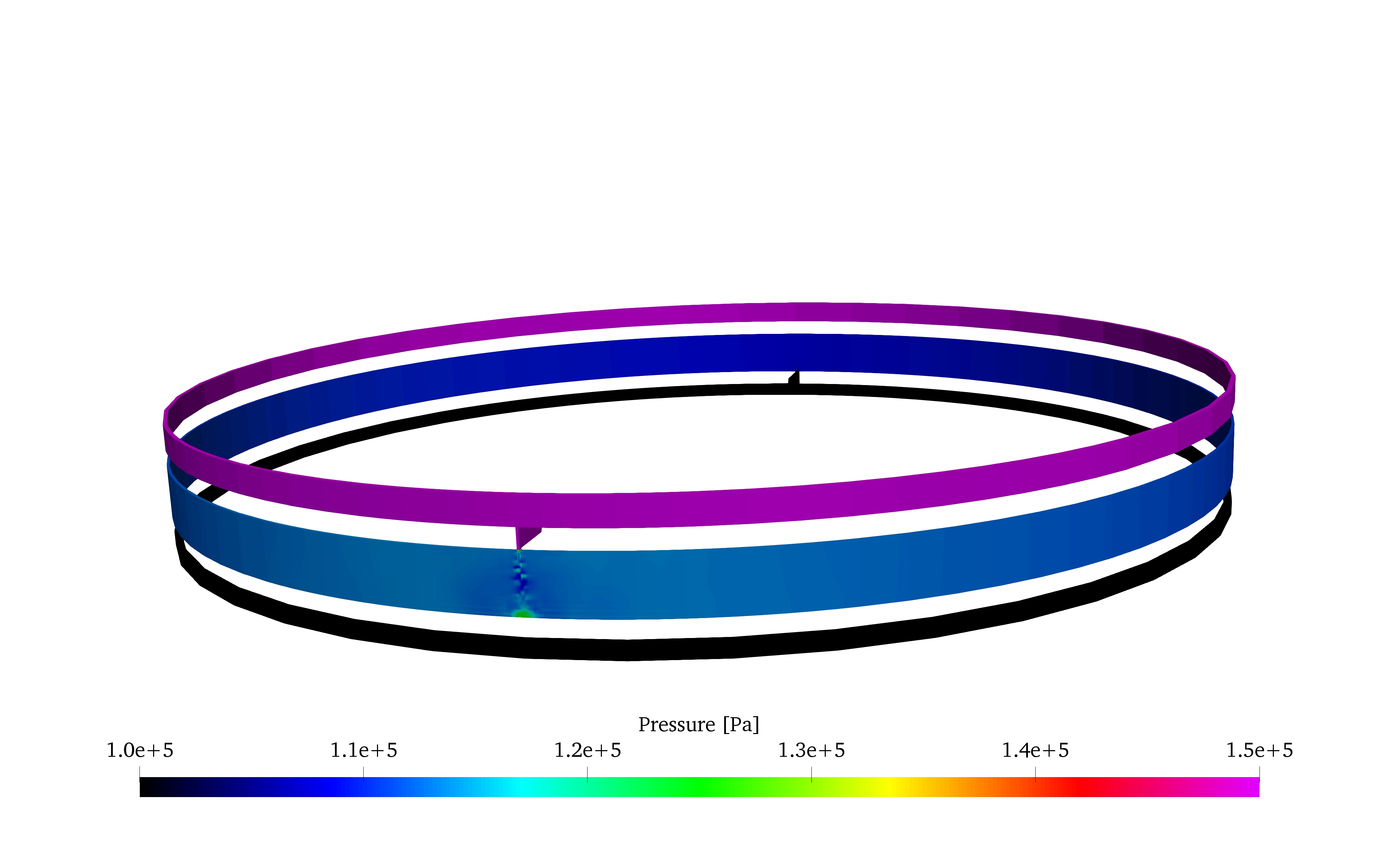} \\
\subfloat[100 time steps FST \label{fig:RingPC}]{\includegraphics[width=0.31\textwidth,trim={0cm 0cm 0cm 0cm},clip]{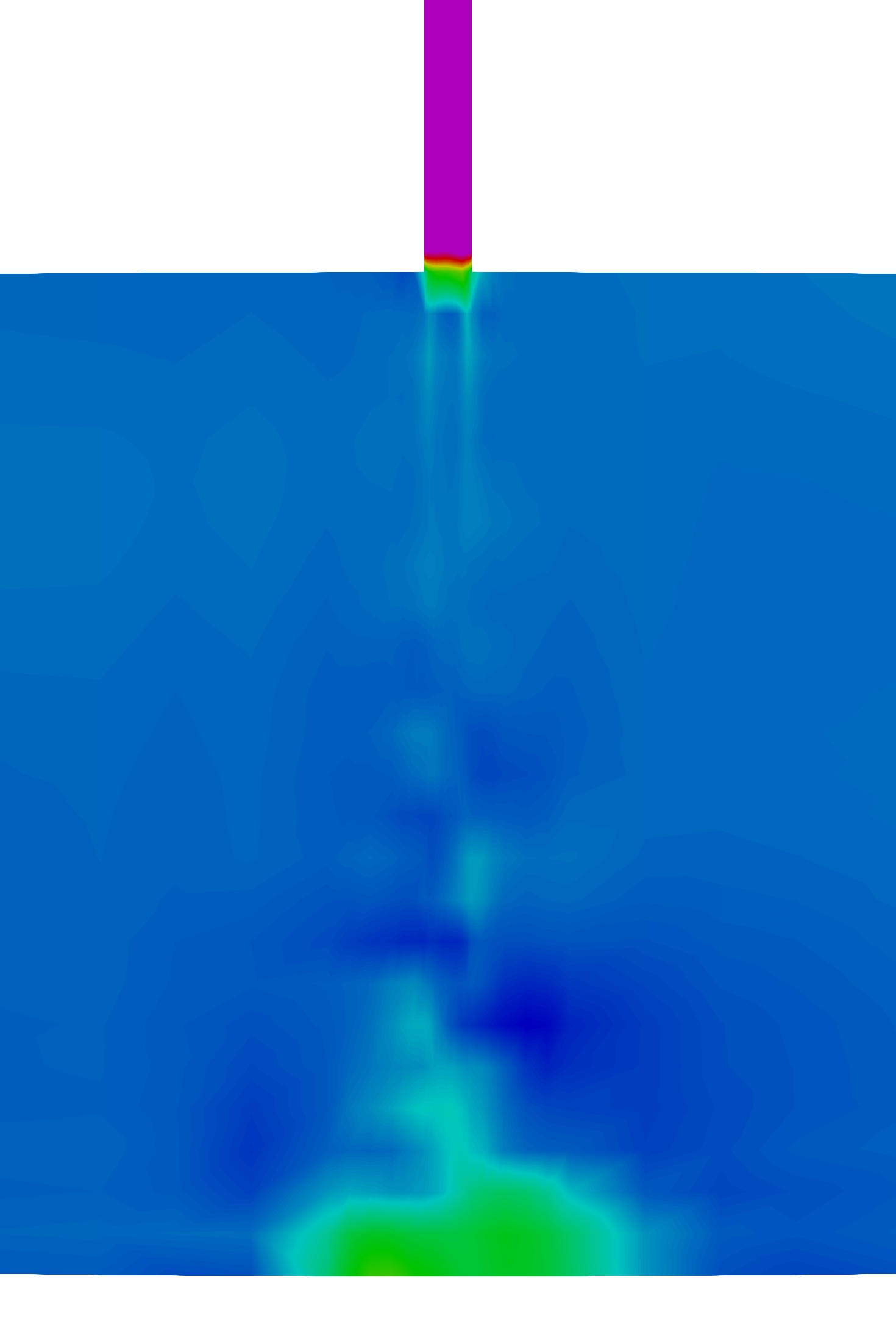}}\,
\subfloat[100 time steps SST \label{fig:RingPM}]{\includegraphics[width=0.31\textwidth,trim={0cm 0cm 0cm 0cm},clip]{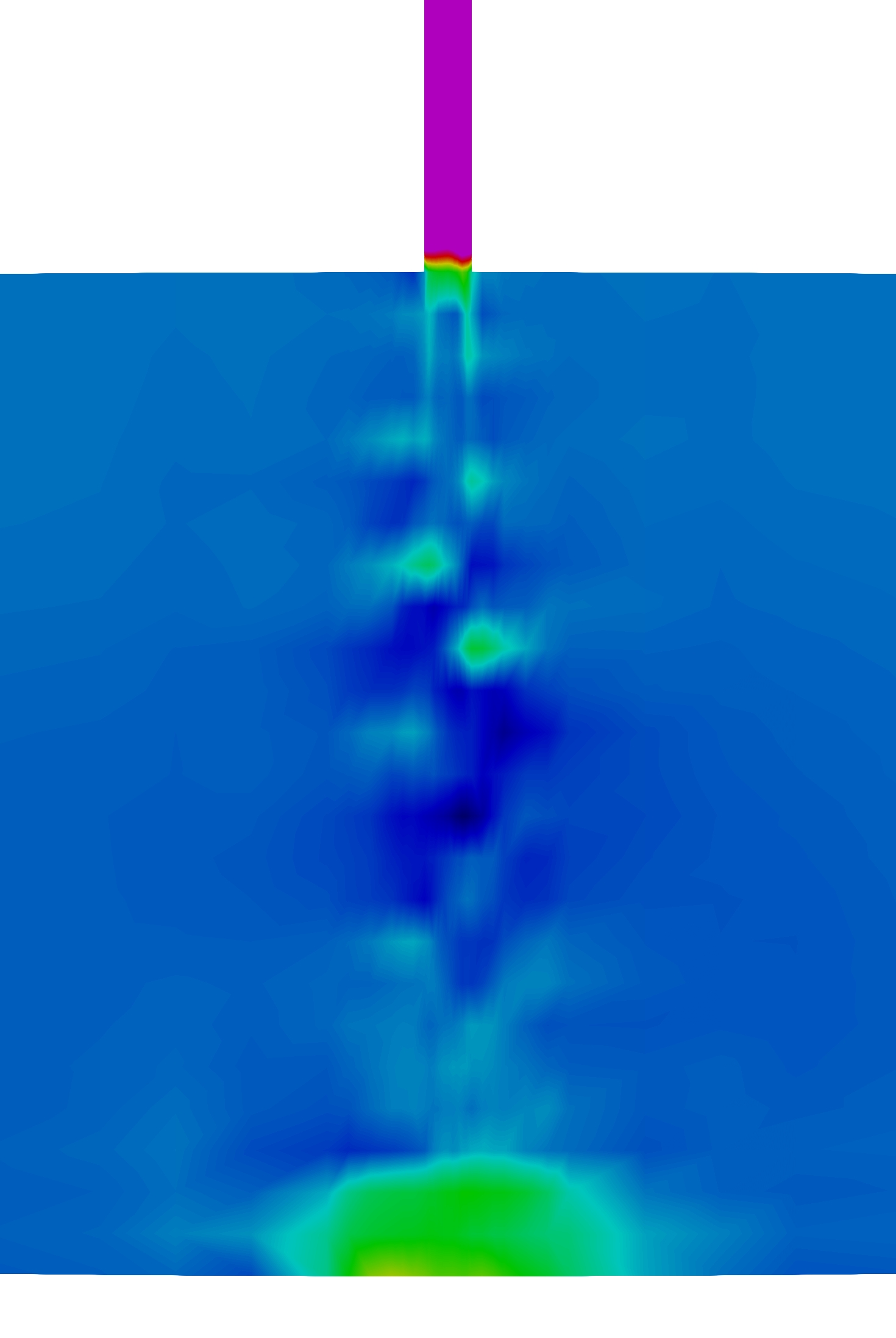}}\,
\subfloat[400 time steps FST \label{fig:RingPF}]{\includegraphics[width=0.31\textwidth,trim={0cm 0cm 0cm 0cm},clip]{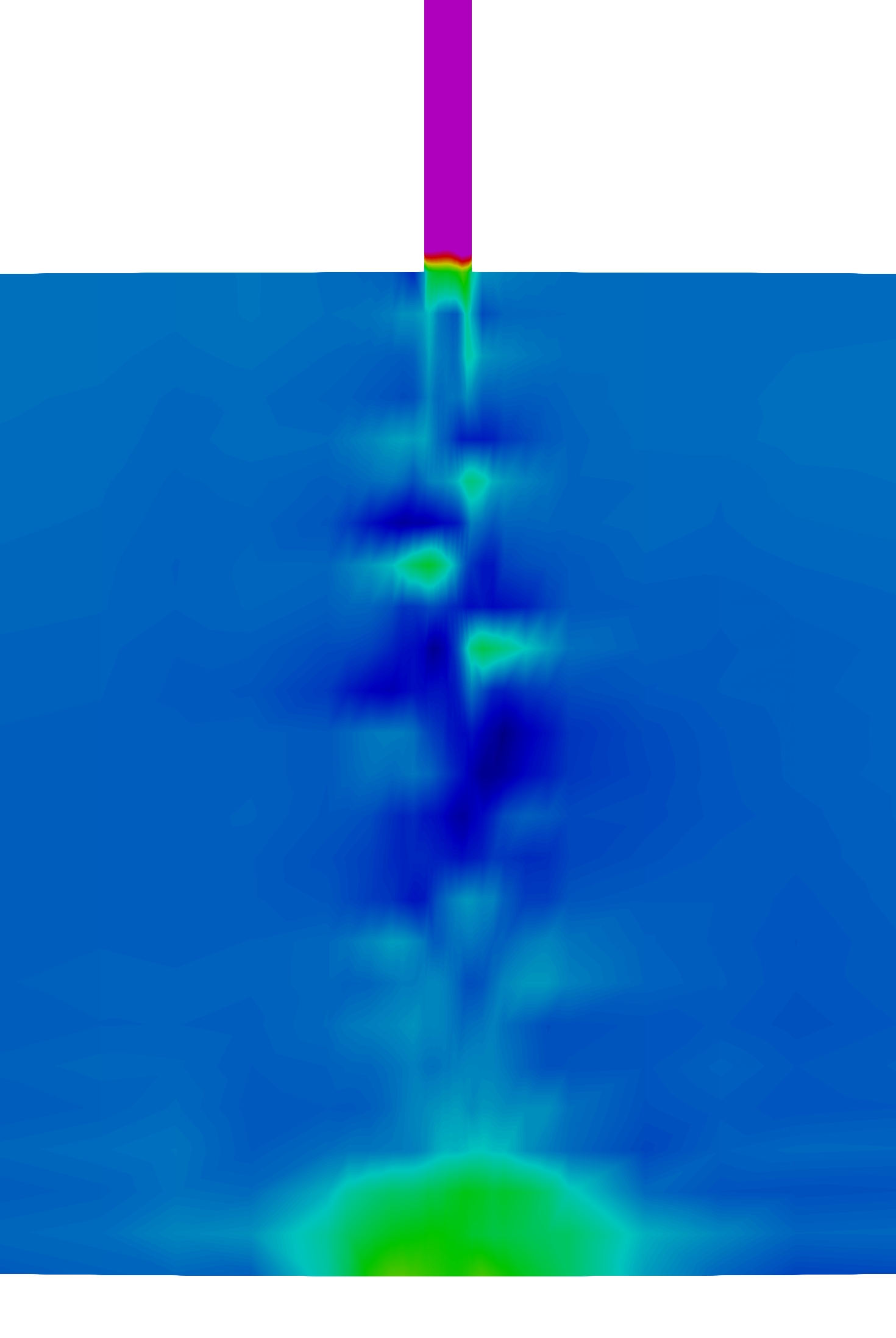}}\,
\caption{Pressure distribution on cylinder wall below upper ring end gap after one degree crank angle for the three different temporal discretizations on the identical spatial discretization.}
\label{fig:RingP}
\end{figure}

Figure~\ref{fig:RingP} compares the pressure distribution on the cylinder wall below the upper ring's end gap for the three different temporal discretizations. Despite the identical spatial discretization, flow features along the upper half of the piston land are resolved by the SST and fine FST discretization, but missing in case of the coarse FST discretization. In that sense, a sufficient temporal accuracy can be obtained with 400 time steps on the FST discretization or 100 time steps on the SST discretization. The four-dimensional time refinement allows us to reduce the number of time steps per degree crank angle ($\nts$), offers an opportunity to save computational time, and motivates a closer look at the simulations' timings.

\begin{table}
\centering
\caption{Mesh data and timings for piston ring blow-by simulations.}
\begin{tabular}{c c c c r r r}
\toprule
 		& $\nts$	  & $\nel$ [Mio] &  $\ndf$[k] & $t_A$ [s] & $t_S$ [s] & $t_T$ [s] \\
  \midrule
Coarse  FST  & 100 & 0.45  &  683 & 362 & 540 & 946\\
SST                  & 100 & 2.91  & 894 & 812 & 1003 & 1898  \\
Fine FST		     & 400 & 0.45 & 683 & 1456 & 2175 & 3761 \\
 \bottomrule
\label{tab:fall}
\end{tabular}
\end{table}

All three simulations have been run in parallel on 120 cores of the Intel Xeon based RWTH cluster using an MPI parallelization. Table~\ref{tab:fall} summarizes the mesh data, as well as timings for the three simulations. Comparing the number of the extruded tetrahedral elements ($\nel$) in the FST meshes to the number of pentatope elements in the SST mesh, one observes a drastic increase. The number of degrees of freedom ($\ndf$), in contrast, is increased only moderately by four-dimensional time refinement. The two main contributions to the total computation time $t_T$, are $t_A$, the time to assemble the linear equation system, and $t_S$, the time to solve the equation system. In sum, we observe a doubling of the computation time $t_T$ comparing the SST simulation to the coarse FST case. Comparing the SST simulation to the fine FST simulation, the computation time is reduced by 50\%. This test case demonstrates how temporally refined simplex space-time finite element meshes can help to reduce the overall simulation time.

\section{Conclusions}
\label{sec:conclusion}
In this paper, we adopted an SUPG stabilized finite element formulation for compressible flows\cite{xu2017compressible} to space-time finite elements. By incorporating a regular simplex element in the definition of the metric tensor $\bG$, we obtained a stabilization matrix $\btau$ that is invariant under permutations of the node numbering of the finite elements. 

We performed compressible flow simulations on several simplex space-time meshes. Our formulation was validated with super- and subsonic test cases with Reynolds numbers ranging from 10 to 1000. SST and UST discretizations proved to be of comparable temporal accuracy as the well-established FST method.

On a valve test case, we demonstrated the capability of the UST method to simulate compressible flows on spatial computational domains undergoing topological changes. Four-dimensional temporal refinement in a pentatope discretization was successfully applied to the transient simulation of blow-by past piston rings of an internal combustion engine.

\section{Acknowledgment}
\label{sec:acknowledgment}
This work was supported by the German Federal Ministry for Economic Affairs and Energy through the Central Innovation Programme for Small and Medium Enterprises in the cooperation project ''Simulationstechnik Kolben-Kolbenring-Zylinder'' under grant number KF3462301PO4. Computing resources were provided by JARA--J\"{u}lich Aachen Research Alliance and RWTH Aachen University IT Center. Furthermore, we want to thank Philipp Knechtges for his helpful remarks during the preparation of the manuscript.


\begin{thebibliography}{10}
\providecommand \doibase [0]{http://dx.doi.org/}%

\bibitem{brooks1982streamline}
Brooks AN, Hughes TJR. Streamline upwind/Petrov-Galerkin formulations for
  convection dominated flows with particular emphasis on the incompressible
  Navier-Stokes equations. {\it Computer Methods in Applied Mechanics and
  Engineering} 1982\string; 32(1-3)\string: 199--259.

\bibitem{hughes1984finite}
Hughes TJR, Tezduyar TE. Finite element methods for first-order hyperbolic
  systems with particular emphasis on the compressible Euler equations. {\it
  Computer methods in applied mechanics and engineering} 1984\string;
  45(1-3)\string: 217--284.

\bibitem{shakib1989finite}
Shakib F. {\it Finite element analysis of the compressible Euler and
  Navier-Stokes equations}. PhD thesis. Stanford University,  1989.

\bibitem{hauke1998comparative}
Hauke G, Hughes TJR. A comparative study of different sets of variables for
  solving compressible and incompressible flows. {\it Computer Methods in
  Applied Mechanics and Engineering} 1998\string; 153(1)\string: 1--44.

\bibitem{xu2017compressible}
Xu F, Moutsanidis G, Kamensky D, et al. Compressible flows on moving domains:
  Stabilized methods, weakly enforced essential boundary conditions, sliding
  interfaces, and application to gas-turbine modeling. {\it Computers \&
  Fluids} 2017\string; 158\string: 201-220.

\bibitem{hughes1988space}
Hughes TJR, Hulbert GM. Space-time finite element methods for elastodynamics:
  formulations and error estimates. {\it Computer Methods in Applied Mechanics
  and Engineering} 1988\string; 66(3)\string: 339--363.

\bibitem{maubach1991iterative}
Maubach JML. {\it Iterative methods for non-linear partial differential
  equations}. PhD thesis. University of Nijmegen,  1991.

\bibitem{idesman2001finite}
Idesman A, Niekamp R, Stein E. Finite elements in space and time for
  generalized viscoelastic maxwell model. {\it Computational Mechanics}
  2001\string; 27(1)\string: 49--60.

\bibitem{behr2008simplex}
Behr M. Simplex space-time meshes in finite element simulations. {\it
  International Journal for Numerical Methods in Fluids} 2008\string;
  57(9)\string: 1421--1434.

\bibitem{neumuller2011refinement}
Neum{\"u}ller M, Steinbach O. Refinement of flexible space--time finite element
  meshes and discontinuous Galerkin methods. {\it Computing and visualization
  in science} 2011\string; 14(5)\string: 189--205.

\bibitem{wang2015discontinuous}
Wang L. {\it Discontinuous galerkin methods on moving domains with large
  deformations}. PhD thesis. University of California, Berkeley,  2015.

\bibitem{wang2015high}
Wang L, Persson PO. A high-order discontinuous Galerkin method with
  unstructured space--time meshes for two-dimensional compressible flows on
  domains with large deformations. {\it Computers \& Fluids} 2015\string;
  118\string: 53--68.

\bibitem{Lehrenfeld2015On}
Lehrenfeld C. {\it On a Space-Time Extended Finite Element Method for the
  Solution of a Class of Two-Phase Mass Transport Problems}. PhD thesis. RWTH
  Aachen University,  2015.

\bibitem{karyofylli2018simplex}
Karyofylli V, Frings M, Elgeti S, Behr M. Simplex space-time meshes in
  two-phase flow simulations. {\it International Journal for Numerical Methods
  in Fluids} 2018\string; 86\string: 218-230.

\bibitem{steinbach2018comparison}
Steinbach O, Yang H. Comparison of algebraic multigrid methods for an adaptive
  space--time finite-element discretization of the heat equation in 3D and 4D.
  {\it Numerical Linear Algebra with Applications} 2018\string; 25(3)\string:
  e2143.
\newblock https://doi.org/10.1002/nla.2143.

\bibitem{elman2014finite}
Elman HC, Silvester DJ, Wathen AJ. {\it Finite elements and fast iterative
  solvers: with applications in incompressible fluid dynamics}.
\newblock Oxford University Press.
\newblock 2005.

\bibitem{pauli2017stabilized}
Pauli L, Behr M. On stabilized space-time FEM for anisotropic meshes:
  Incompressible Navier--Stokes equations and applications to blood flow in
  medical devices. {\it International Journal for Numerical Methods in Fluids}
  2017\string; 85\string: 189--209.

\bibitem{buchholz1992perfect}
Buchholz RH. Perfect pyramids. {\it Bulletin of the Australian Mathematical
  Society} 1992\string; 45(3)\string: 353--368.

\bibitem{hauke2001simple}
Hauke G. Simple stabilizing matrices for the computation of compressible flows
  in primitive variables. {\it Computer Methods in Applied Mechanics and
  Engineering} 2001\string; 190(51)\string: 6881--6893.

\bibitem{knechtges2018simulation}
Knechtges P. {\it Simulation of Viscoelastic Free-Surface Flows}. PhD thesis.
  RWTH Aachen University,  2018.

\bibitem{higham1987computing}
Higham NJ. Computing real square roots of a real matrix. {\it Linear Algebra
  and its applications} 1987\string; 88\string: 405--430.

\bibitem{nagLibrary}
{The NAG Library}. The Numerical Algorithms Group (NAG), Oxford, United
  Kingdom. www.nag.com

\bibitem{behr1994finite}
Behr M, Tezduyar TE. Finite element solution strategies for large-scale flow
  simulations. {\it Computer Methods in Applied Mechanics and Engineering}
  1994\string; 112(1-4)\string: 3--24.

\bibitem{anderson2010fundamentals}
Anderson~Jr JD. {\it Fundamentals of aerodynamics}.
\newblock Tata McGraw-Hill Education.
\newblock 2010.

\bibitem{carter1972numerical}
Carter JE. Numerical solutions of the Navier-Stokes equations for the
  supersonic laminar flow over a two-dimensional compression corner. tech.
  rep., NASA;  1972.

\bibitem{thompson1971compressible}
Thompson PA. {\it Compressible-fluid dynamics}.
\newblock McGraw-Hill.
\newblock 1971.

\bibitem{oliva2016numerical}
Oliva A, Held S. Numerical multiphase simulation and validation of the flow in
  the piston ring pack of an internal combustion engine. {\it Tribology
  International} 2016\string; 101\string: 98--109.

\bibitem{geuzaine2009gmsh}
Geuzaine C, Remacle JF. Gmsh: A 3-D finite element mesh generator with built-in
  pre-and post-processing facilities. {\it International Journal for Numerical
  Methods in Engineering} 2009\string; 79(11)\string: 1309--1331.

\end{thebibliography}
\end{document}